\documentclass[12pt]{article}
\usepackage{graphicx}
\usepackage{amssymb}
\usepackage{latexsym}

\def\jl#1{\bigskip{\bf #1}\bgroup\it\ }
\def\ejl{\egroup\par\bigskip}

\begin{document}

\title{\Large \bf
Detection of some elements in the stable homotopy groups of spheres
\thanks{\footnotesize
This work was supported partially by the National Natural Science
Foundation of China (Grant No. 10501045), and the Fund of the
Personnel Division of Nankai University.}} {
\author {\begin{tabular}{c}
\sc Liu Xiugui\\
\end{tabular} \\
{School of Mathematical Sciences and LPMC}\\
{Nankai University}\\
{Tianjin 300071, P. R. China}\\
\\
{E-mail:\quad xgliu@nankai.edu.cn} }\date{} \maketitle }

\baselineskip=16pt
\begin{center}
\begin{minipage}{135mm}
\noindent {\bf Abstract} In this paper we constructs a new
nontrivial family in the stable homotopy groups
of spheres $\pi_{p^nq+2pq+q-3}S$ which is of order $p$ and is
represented by $k_0h_{n} \in Ext_A^{3,p^nq+2pq+q}(\mathbb{Z}_p,\mathbb{Z}_p)$ in the
Adams spectral sequence, where $p\geq 5$ is an odd prime, $n\geq 3$ and
$q=2(p-1)$. In the course of the proof, a new family of homotopy
elements in $\pi_{\ast}V(1)$ which is represented by
$\beta_{\ast}{i^{\prime}}_{\ast}i_{\ast}({h}_n)\in
Ext_A^{2,p^nq+(p+1)q+1}(H^{\ast}V(1),\mathbb{Z}_p)$ in the Adams sequence is
detected.
\\


\hspace{2.3cm}May spectral sequence, Steenrod algebra.

{\bf 2000 MR Subject Classification} 55Q45
\end{minipage}
\end{center}
\vspace {0.5cm}
{\large \bf 1 Introduction and the main results}
\\
\\
\noindent Let $A$ be the mod $p$ Steenrod algebra and $S$ be the
sphere spectrum localized at an odd prime $p$. To determine the
stable homotopy groups of spheres $\pi_{\ast}S$ is one of the
central problems in homotopy theory. One of the main tools to reach
it is the Adams spectral sequence:
$$E_2^{s,t}=Ext_A^{s,t}(\mathbb{Z}_p , \mathbb{Z}_p)\Rightarrow
\pi_{t-s}S,$$ where the $E_2^{s,t}$-term is the cohomology of $A$.
So far, not so many families of homotopy elements in $\pi_{\ast}S$
have been detected. For example, a family
$\varsigma_{n-1}\in\pi_{p^nq+q-3}S$ for $n\geq 2$ which has
filtration $3$ in the Adams spectral sequence and is represented by
$h_0b_{n-1}\in Ext_A^{3,p^nq+q}(\mathbb{Z}_p , \mathbb{Z}_p)$ has
been detected in reference [1], where $q=2(p-1)$. In this paper, we
also detect a family of homotopy elements in $\pi_{p^nq+pq-3}S$
which has filtration $3$ and is represented by $k_0h_n\in
Ext_A^{3,p^nq+2pq+q}(\mathbb{Z}_p , \mathbb{Z}_p)$ in the Adams
spectral sequence.

From reference [2], $Ext_A^{1,\ast}(\mathbb{Z}_p,\mathbb{Z}_p)$ has $\mathbb{Z}_p$-bases consisting of
$a_0\in Ext_A^{1,1}(\mathbb{Z}_p , \mathbb{Z}_p)$, $h_i\in Ext_A^{1,p^iq}(\mathbb{Z}_p , \mathbb{Z}_p)$
for all $i\geq 0$ and $Ext_A^{2,\ast}(\mathbb{Z}_p , \mathbb{Z}_p)$ has $\mathbb{Z}_p$-bases
consisting of $\alpha_2$, $a_0^2$, $a_0h_i(i>0)$, $g_i(i\geq 0)$, $k_i(i\geq 0)$,
$e_i(i\geq 0)$, and $h_ih_j(j\geq i+2,i\geq 0)$ whose internal degrees are $2q+1$, $2$,
$p^iq+1,p^{i+1}q+2p^iq$, $2p^{i+1}+p^iq$, $p^{i+1}q$ and $p^iq+p^jq$
respectively.

Let $M$ be the Moore spectrum modulo a prime $p\geq 5$ given by
the cofibration
$$S\stackrel{p}{\rightarrow}S\stackrel{i}{\rightarrow}M\stackrel{j}{\rightarrow}\Sigma S.
\eqno{(1.1)}$$
Let $\alpha:\Sigma^q M\rightarrow M$ be the Adams map and $K$ be its
cofibre given by the cofibration
$$\Sigma^qM\stackrel{\alpha}{\rightarrow}M\stackrel{i^{\prime}}{\rightarrow}
K\stackrel{j^{\prime}}{\rightarrow}\Sigma^{q+1}M
,\eqno{(1.2)}$$
where $q=2(p-1)$. This spectrum which we briefly write as $K$ is
known to be the Toda-Smith spectrum $V(1)$. Let $V(2)$ be the
cofibre of $\beta:\Sigma^{(p+1)q}K\rightarrow K$ given by the
cofibration

$$\Sigma^{(p+1)q}K\stackrel{\beta}{\rightarrow}K\stackrel
{\overline{i}}{\rightarrow}V(2)\stackrel{\overline{j}}{\rightarrow}\Sigma^{(p+1)q+1}K.
\eqno{(1.3)}
$$

Our results can be stated as follows.

{\bf Theorem I} {\it Let $p\geq 5,n\geq 3$,then
$$\beta_{\ast}{i^{\prime}}_{\ast}i_{\ast}({h}_n)\not=0\in
Ext_A^{2,p^nq+(p+1)q+1}(H^{\ast}K,\mathbb{Z}_p)$$
is a permanent cycle in the Adams spectral sequence and converges to a
nontrivial element $\zeta_n \in \pi_{p^nq+(p+1)q-1}K$, where $h_n\in
Ext_A^{1,p^nq}(\mathbb{Z}_p,\mathbb{Z}_p)$.}

{\bf Theorem II} {\it Let $p\geq 5$, $n\geq 3$, then
$$k_0h_n\not=0\in Ext_A^{3,p^nq+2pq+q}(\mathbb{Z}_p ,
\mathbb{Z}_p)$$ is a permanent cycle in the Adams spectral sequence and it converges to a
nontrivial element of order $p$ in $\pi_{p^nq+2pq+q-3}S$.}

{\bf Remark:} The $k_0h_n$-element obtained in Theorem II is an indecomposable element in $\pi_{\ast}S$,
i.e., it is not a composition of elements of lower filtration in $\pi_{\ast}S$, because $h_n$($n>0$)
is known to die in the Adams spectral sequence.

After giving some useful propositions in Section 2, the proofs of
the main theorems will be given in Section 3.
\\
\\
{\large\bf 2 Some preliminaries on low-dimensional Ext groups}
\\
\\
\noindent In this section, we will prove some results on Ext groups
of lower dimension which will be used in the proofs of the theorems.

{\bf Proposition 2.1} {\it Let $p\geq 5,n\geq 3, a_0\in Ext_A^{1,1}(\mathbb{Z}_p,\mathbb{Z}_p)$,$h_n\in Ext_A^{1,p^nq}(\mathbb{Z}_p,\mathbb{Z}_p)$,
$b_n\in Ext_A^{2,p^{n+1}q}(\mathbb{Z}_p,\mathbb{Z}_p)$ respectively. Then we have the following:

(1) $Ext_A^{4,p^nq+pq+2}(\mathbb{Z}_p,\mathbb{Z}_p)\cong \mathbb{Z}_p\{a_0a_0h_1h_n\}$.

(2) $Ext_A^{5,p^nq+(p+2)q}(\mathbb{Z}_p,\mathbb{Z}_p)=0$.}

{\it Proof.} (1) See [3, Theorem 4.1].

(2) The proof is similar to that given in the proof of [4,
Proposition 1.2]. We can show that in the May spectral sequence
$E_1^{5, p^nq+(p+2)q,\ast}=0$. Then
$$Ext_A^{5,p^nq+(p+2)q}(\mathbb{Z}_p,\mathbb{Z}_p)=0.$$ Here the proof is omitted.

{\protect\vspace{-15pt}\rightline{$\square$}}

The following lemma is used in the proofs of many propositions in
this section.

First recall spectra $V(k)=\{V(k)_n\}$ for $n\geq -1$ which are
so-called Toda-Smith spectra. The spectrum $V(n)$ given in [5] such
that the $Z_p$-cohomology
$$H^{\ast}(V(n),\mathbb{Z}_p)\cong E(n)=E(
\mathcal{Q}_0, \mathcal{Q}_1, \cdots, \mathcal{Q}_n ),$$ the
exterior algebra generator by Milnor basis elements $\mathcal{Q}_0$,
$\mathcal{Q}_1$, $\cdots$, $\mathcal{Q}_n$ in $A$. The spectra
$V(n)$ for $n\geq -1$ are defined inductively by $V(-1)=S$ and the
cofibration
$$
\Sigma^{2(p^n-1)}V(n-1)\stackrel{\alpha^{(n)}}{\rightarrow}V(n-1)
\stackrel{i_{n}}{\rightarrow}V(n)
\stackrel{j_n}{\rightarrow}\Sigma^{2p^n-1}V(n-1). \eqno{(2.1)}$$
When $n=0$, $1$, $2$, the above cofibration sequence just is the
cofibration sequences (1.1), (1.2) and (1.3) respectively.
$\alpha^{(n)}$ stand for the maps $p$, $\alpha$, $\beta$ in (1.1),
(1.2) and (1.3) respectively. Here $V(-1)=S, V(0)=M, V(1)=K, i_0=i,
i_1=i^{\prime}, i_2=\bar{i}, j_0=j, j_1=j^{\prime}, j_2=\bar{j}$.
The existence of $V(n)$ is assured [5, Theorem 1.1] for $n=1$,
$p\geq 3$ and for $n=2$, $p\geq 5$.

By the definition of $Ext$ groups, from (2.1) we can easily have the
following lemma.

{\bf Lemma 2.1} {\it With notations as above. We have the following
two long exact sequence:

(1) $\cdots\rightarrow Ext_A^{s-1, t-(2p^n-1)}(H^{\ast}V(n-1),\Box)
\stackrel{\alpha^{(n)}_{\ast}}{\rightarrow}
Ext_A^{s,t}(H^{\ast}V(n-1),\Box)
\stackrel{(i_n)_{\ast}}{\rightarrow}
\\Ext_A^{s,t}(H^{\ast}V(n),\Box)
\stackrel{(j_n)_{\ast}}{\rightarrow}
Ext_A^{s,t-(2p^n-1)}(H^{\ast}V(n),\Box) \rightarrow\cdots $.

(2) $\cdots\rightarrow Ext_A^{s-1, t-(2p^n-1)}(\Box, H^{\ast}V(n-1))
\stackrel{\alpha^{(n)^\ast}}{\rightarrow}
Ext_A^{s,t}(\Box,H^{\ast}V(n-1))
\stackrel{(j_n)^{\ast}}{\rightarrow} Ext_A^{s,t}(\Box,H^{\ast}V(n))
\stackrel{(i_n)^{\ast}}{\rightarrow} Ext_A^{s,t}(\Box,H^{\ast}V(n))
\rightarrow\cdots $.

Here $\Box$ are an arbitrary $A$-module.}

{\protect\vspace{-15pt}\rightline{$\square$}}

{\bf Proposition 2.2} {\it Let $p\geq 5$, $n\geq 3$. Then
$Ext_A^{3,p^nq+(p+2)q+1}(H^{\ast}M,H^{\ast}M)$ has a unique
generator
 $\overline{\overline{h_ng_0}}$, where $\overline{\overline{h_ng_0}}$ satisfies
 $i^{\ast}j_{\ast}\overline{\overline{h_ng_0}}=h_ng_0$, the generator of
$Ext_A^{3,p^nq+(p+2)q}(\mathbb{Z}_p,\mathbb{Z}_p)$ stated in [6,
Table 8.1].}

{\it Proof.} First consider the exact sequence
$$
\begin{array}{r}Ext_A^{3,p^nq+(p+2)q+1}(\mathbb{Z}_p,\mathbb{Z}_p)\stackrel{j^{\ast}}{\rightarrow}
Ext_A^{3,p^nq+(p+2)q}(\mathbb{Z}_p,H^{\ast}M)\stackrel{i^{\ast}}{\rightarrow}\\
Ext_A^{3,p^nq+(p+2)q}(\mathbb{Z}_p,\mathbb{Z}_p)\stackrel{p^{\ast}}{\rightarrow}
Ext_A^{4,p^nq+(p+2)q+1}(\mathbb{Z}_p,\mathbb{Z}_p)
\end{array}
$$
induced by (1.1). Since we know that
$Ext_A^{3,p^nq+(p+2)q+1}(\mathbb{Z}_p,\mathbb{Z}_p)$ is zero (cf.
[6, Table 8.1]) and
$Ext_A^{4,p^nq+(p+2)q+1}(\mathbb{Z}_p,\mathbb{Z}_p)$ is zero (cf.
[7, Proposition 2.1]), the above $i^{\ast}$ is an isomorphism. Then
we see that $Ext_A^{3,p^nq+(p+2)q}(\mathbb{Z}_p,H^{\ast}M)$ has
unique generator $\overline{h_ng_0}$, where $\overline{h_ng_0}$
satisfies $i^{\ast}\overline{h_ng_0}=h_ng_0$, the unique generator
of $Ext_A^{3,p^nq+(p+2)q}(\mathbb{Z}_p,\mathbb{Z}_p)$ stated in [6,
Table 8.1].

At last, look at the following exact sequence induced by (1.1)
$$\begin{array}{r}
Ext_A^{3,p^nq+(p+2)q+1}(\mathbb{Z}_p,H^{\ast}M)\stackrel{i_{\ast}}{\rightarrow}
Ext_A^{3,p^nq+(p+2)q+1}(H^{\ast}M,H^{\ast}M)
\stackrel{j_{\ast}}{\rightarrow}\\
Ext_A^{3,p^nq+(p+2)q}(\mathbb{Z}_p,H^{\ast}M)
\stackrel{p_{\ast}}{\rightarrow}Ext_A^{4,p^nq+(p+2)q+1}(\mathbb{Z}_p,H^{\ast}M).
\end{array}
$$
Since the first group is zero by virtue of
$Ext_A^{3,p^nq+(p+2)q+r}(\mathbb{Z}_p,\mathbb{Z}_p)=0$ for $r=1,2$
(cf. [6, Table 8.1]) and the forth group is zero by virtue of the
facts that $Ext_A^{4,p^nq+(p+2)q+t}(\mathbb{Z}_p,\mathbb{Z}_p)=0$
for $t=1,2$ (cf. [7, Proposition 2.1]), then the above $j_{\ast}$ is
an isomorphism. Thus $Ext_A^{3,p^nq+(p+2)q+1}(H^{\ast}M,H^{\ast}M)$
has a unique generator $\overline{\overline{h_ng_0}}$, where
$\overline{\overline{h_ng_0}}$ satisfies
$j_{\ast}\overline{\overline{h_ng_0}}=\overline{h_ng_0}$. This
finishes the proof of Proposition 2.2.

{\protect\vspace{-15pt}\rightline{$\square$}}

{\bf Proposition 2.3} {\it Let $p\geq 5,n\geq 3$, then
$$Ext_A^{3,p^nq+(p+2)q}(H^{\ast}M,H^{\ast}M)
\cong \mathbb{Z}_p\{i_{\ast}j_{\ast}\overline{\overline{h_ng_0}},
j^{\ast}i^{\ast}\overline{\overline{h_ng_0}}\}.$$}
\indent{\it Proof.} Consider the exact sequence
$$\begin{array}{r}
Ext_A^{2,p^nq+(p+2)q-1}(\mathbb{Z}_p,\mathbb{Z}_p)\stackrel{p^{\ast}}{\rightarrow}
Ext_A^{3,p^nq+(p+2)q}(\mathbb{Z}_p,\mathbb{Z}_p)
\stackrel{j^{\ast}}{\rightarrow}\\
Ext_A^{3,p^nq+(p+2)q-1}(\mathbb{Z}_p,H^{\ast}M)
\stackrel{i^{\ast}}{\rightarrow}
Ext_A^{3,p^nq+(p+2)q-1}(\mathbb{Z}_p,\mathbb{Z}_p)
\end{array}
$$
induced by (1.1). Since
$Ext_A^{2,p^nq+(p+2)q-1}(\mathbb{Z}_p,\mathbb{Z}_p)=0
=Ext_A^{3,p^nq+(p+2)q-1}(\mathbb{Z}_p,\mathbb{Z}_p)$ (cf. [6, Table
8.1]), the above $j^{\ast}$ is an isomorphism. Moreover we also know
$Ext_A^{3,p^nq+(p+2)q}(\mathbb{Z}_p,\mathbb{Z}_p)=\mathbb{Z}_p\{h_ng_0\}$
(cf. [6, Table 8.1]). Thus we can have that
$Ext_A^{3,p^nq+(p+2)q-1}(\mathbb{Z}_p,H^{\ast}M)=\mathbb{Z}_p\{j^{\ast}(h_ng_0)\}$.

Now observe the following exact sequence
$$\begin{array}{r}
Ext_A^{2,p^nq+(p+2)q-1}(\mathbb{Z}_p,H^{\ast}M)\stackrel{p_{\ast}}{\rightarrow}
Ext_A^{3,p^nq+(p+2)q}(\mathbb{Z}_p,H^{\ast}M)
\stackrel{i_{\ast}}{\rightarrow}\\
Ext_A^{3,p^nq+(p+2)q}(H^{\ast}M,H^{\ast}M)
\stackrel{j_{\ast}}{\rightarrow}
Ext_A^{3,p^nq+(p+2)q-1}(\mathbb{Z}_p,H^{\ast}M)\stackrel{p_{\ast}}{\rightarrow}
\end{array}$$
induced by (1.1). Since
$Ext_A^{2,p^nq+(p+2)q+r}(\mathbb{Z}_p,\mathbb{Z}_p)=0$ for $r=-1,0$
(cf. [2]), we can easily get that
$Ext_A^{2,p^nq+(p+2)q-1}(\mathbb{Z}_p,H^{\ast}M)=0$. By virtue of
the fact
$Ext_A^{3,p^nq+(p+2)q-1}(\mathbb{Z}_p,H^{\ast}M)=\mathbb{Z}_p\{j^{\ast}(h_ng_0)\}$,
we have that the image of the second $p_{\ast}$ is
$p_{\ast}j^{\ast}(h_ng_0)=j^{\ast}p_{\ast}(h_ng_0)=j^{\ast}p^{\ast}(g_0h_n)=0$.

From the facts that $$Ext_A^{3,p^nq+(p+2)q}(\mathbb{Z}_p,H^{\ast}M)\cong
\mathbb{Z}_p\{\overline{h_ng_0}\}\cong \mathbb{Z}_p\{j_{\ast}\overline{\overline{h_ng_0}}\}
$$
and $$Ext_A^{3,p^nq+(p+2)q-1}(\mathbb{Z}_p,H^{\ast}M)\cong \mathbb{Z}_p\{j^{\ast}(h_ng_0)\}\cong
\mathbb{Z}_p\{j^{\ast}i^{\ast}j_{\ast}
\overline{\overline{h_ng_0}}\}\cong
\mathbb{Z}_p\{j_{\ast}j^{\ast}i^{\ast}\overline{\overline{h_ng_0}}\},$$
we can easily get that $Ext_A^{3,p^nq+(p+2)q}(H^{\ast}M,H^{\ast}M)
\cong \mathbb{Z}_p\{i_{\ast}j_{\ast}\overline{\overline{h_ng_0}},
j^{\ast}i^{\ast}\overline{\overline{h_ng_0}}\}$. This shows Proposition
2.3.

{\protect\vspace{-15pt}\rightline{$\square$}}

{\bf Proposition 2.4} {\it Let $p\geq 5,n\geq 3$, then we have

(1) $i^{\ast}d_2(i_{\ast}j_{\ast}(\overline{\overline{h_ng_0}}))\not= 0$.

(2) $d_{2}(j^{\ast}i^{\ast}(\overline{\overline{h_ng_0}}))\not=0$, where
$$d_2: Ext_A^{3,p^nq+(p+2)q}(H^{\ast}M,H^{\ast}M)\rightarrow
Ext_A^{5,p^nq+(p+2)q+1}(H^{\ast}M,H^{\ast}M)$$ is the differential of the Adams spectral sequence.}

{\it Proof.} (1) From [7, p.488] we know that
$d_2(i_{\ast}(h_ng_0))\not=0$. By Proposition 2.2,
$d_2(i_{\ast}(h_ng_0))=d_2(i_{\ast}i^{\ast}j_{\ast}(\overline{\overline{h_ng_0}}))=
d_2(i^{\ast}i_{\ast}j_{\ast}(\overline{\overline{h_ng_0}}))=
i^{\ast}d_2(i_{\ast}j_{\ast}(\overline{\overline{h_ng_0}}))$. The
desired result follows.

(2) Consider the exact sequence
$$\begin{array}{r}
Ext_A^{4,p^nq+(p+2)q}(\mathbb{Z}_p,\mathbb{Z}_p)\stackrel{p^{\ast}}{\rightarrow}
Ext_A^{5,p^nq+(p+2)q+1}(\mathbb{Z}_p,\mathbb{Z}_p)
\stackrel{j^{\ast}}{\rightarrow}Ext_A^{5,p^nq+(p+2)q}(\mathbb{Z}_p,H^{\ast}M)\\
\stackrel{i^{\ast}}{\rightarrow}
Ext_A^{5,p^nq+(p+2)q}(\mathbb{Z}_p,\mathbb{Z}_p)
\end{array}$$ induced by (1.1). We claim that the above $j^{\ast}$ is an
isomorphism. By virtue of the fact that
$Ext_A^{5,p^nq+(p+2)q}(\mathbb{Z}_p,\mathbb{Z}_p)=0$ (cf.
Proposition 2.1), we see that the above $j^{\ast}$ is an
epimorphism. Note that
$Ext_A^{4,p^nq+(p+2)q}(\mathbb{Z}_p,\mathbb{Z}_p) \cong
\mathbb{Z}_p\{g_0b_{n-1}\}$ (cf. [7, Proposition 2.1]). Since
$p^{\ast}(g_0b_{n-1})=a_0g_0b_{n-1}=0$ (Note: $a_0g_0=0$ by [6,
Table 8.2]), so $kerj^{\ast}=imp^{\ast}=0$, i.e., the above
$j^{\ast}$ is a monomorphism. The proof of the claim is finished.
Since $\alpha_2b_0h_n\not=0\in
Ext_A^{5,p^nq+(p+2)q+1}(\mathbb{Z}_p,\mathbb{Z}_p)$ (cf. [7,
Proposition 2.1]), so by the claim we get that
$j^{\ast}(\alpha_2b_0h_n)\not= 0\in
Ext_A^{5,p^nq+(p+2)q}(\mathbb{Z}_p,H^{\ast}M)$. Note the fact that
$d_2(h_ng_0)=\alpha_2b_0h_n\not=0$, so $j^{\ast}d_2
(j_{\ast}i^{\ast}\overline{\overline{h_ng_0}})=j_{\ast}d_2(j^{\ast}i^{\ast}
\overline{\overline{h_ng_0}}) \not=0$ by Proposition 2.2. Thus
$d_{2}(j^{\ast}i^{\ast}(\overline{\overline{h_ng_0}}))\not=0$.

{\protect\vspace{-15pt}\rightline{$\square$}}

{\bf Proposition 2.5} {\it Let $p\geq 5,n\geq 3$,
then $$Ext_A^{3,p^nq+(p+1)q+2}(H^{\ast}K,H^{\ast}M)=0.$$}
\indent{\it Proof.} Consider the exact sequence
$$\begin{array}{r}
Ext_A^{3,p^nq+(p+1)q+3}(H^{\ast}M,\mathbb{Z}_p)
\stackrel{j^{\ast}}{\rightarrow}Ext_A^{3,p^nq+(p+1)q+2}
(H^{\ast}M,H^{\ast}M)\stackrel{i^{\ast}}{\rightarrow}\\
Ext_A^{3,p^nq+(p+1)q+2}(H^{\ast}M,\mathbb{Z}_p).
\end{array}
$$
Since the first and third group are zero by the facts
$Ext_A^{3,p^nq+(p+1)q+r}(\mathbb{Z}_p,\mathbb{Z}_p)=0$ for $r=1,2,3$
(cf. [6, Table 8.1]), so the second group is zero.

Look at the exact sequence
$$\begin{array}{cl}
Ext_A^{3,p^nq+pq+2}(\mathbb{Z}_p,\mathbb{Z}_p)\stackrel{i_{\ast}}{\rightarrow}Ext_A^{3,p^nq+pq+2}(H^{\ast}M,\mathbb{Z}_p)
&\stackrel{j_{\ast}}{\rightarrow}Ext_A^{3,p^nq+pq+1}(\mathbb{Z}_p,\mathbb{Z}_p)\\
&\stackrel{p_{\ast}}{\rightarrow}
Ext_A^{4,p^nq+pq+2}(\mathbb{Z}_p,\mathbb{Z}_p)
\end{array}$$ induced by (1.1).
Since we know that
$Ext_A^{3,p^nq+pq+1}(\mathbb{Z}_p,\mathbb{Z}_p)\cong
\mathbb{Z}_p\{a_0h_1h_n\}$ (cf. [6, Table 8.1]) and
$Ext_A^{4,p^nq+pq+2}(\mathbb{Z}_p,\mathbb{Z}_p) \cong
\mathbb{Z}_p\{a_0a_0h_1h_n\}$ by Proposition 2.1, so the above
$p_{\ast}$ is an isomorphism. $im j_{\ast}=0$ since $p_{\ast}$ is an
isomorphism. $im i_{\ast}=0$ by the fact that
$Ext_A^{3,p^nq+pq+2}(\mathbb{Z}_p,\mathbb{Z}_p)=0$ (cf. [6, Table
8.1]). Thus we can have that
$Ext_A^{3,p^nq+pq+2}(H^{\ast}M,\mathbb{Z}_p)=0$.

Observe the following exact sequence induced by (1.1)
$$\begin{array}{r}
Ext_A^{2,p^nq+pq}(\mathbb{Z}_p,\mathbb{Z}_p)\stackrel{p_{\ast}}{\rightarrow}Ext_A^{3,p^nq+pq+1}(\mathbb{Z}_p,\mathbb{Z}_p)
\stackrel{i_{\ast}}{\rightarrow}Ext_A^{3,p^nq+pq+1}(H^{\ast}M,\mathbb{Z}_p)\\\stackrel{j_{\ast}}{\rightarrow}
Ext_A^{3,p^nq+pq}(\mathbb{Z}_p,\mathbb{Z}_p)\stackrel{p_{\ast}}{\rightarrow}Ext_A^{4,p^nq+pq+1}(\mathbb{Z}_p,\mathbb{Z}_p).
\end{array}$$
Since $Ext_A^{2,p^nq+pq}(\mathbb{Z}_p,\mathbb{Z}_p)\cong
\mathbb{Z}_p\{h_1h_n\}$ and
$Ext_A^{3,p^nq+pq+1}(\mathbb{Z}_p,\mathbb{Z}_p)\cong
\mathbb{Z}_p\{a_0h_1h_n\}$ (cf. [6, Table 8.1]), we know that the
first $p_{\ast}$ is an isomorphism. Similarly by virtue of the facts
that $Ext_A^{3,p^nq+pq}(\mathbb{Z}_p,\mathbb{Z}_p)\cong
\mathbb{Z}_p\{b_0h_n,h_1b_{n-1}\}$ (cf. [6, Table 8.1]) and
$Ext_A^{4,p^nq+pq+1}(\mathbb{Z}_p,\mathbb{Z}_p)\cong
\mathbb{Z}_p\{a_0 b_0h_n,a_0h_1b_{n-1}\}$ (cf. [7, Proposition
2.1]), we get that the second $p_{\ast}$ is also an isomorphism.
Thus $Ext_A^{3,p^nq+pq+1}(H^{\ast}M,\mathbb{Z}_p)=0$.

Look at the exact sequence
$$\begin{array}{r}
0=Ext_A^{3,p^nq+pq+2}(H^{\ast}M,\mathbb{Z}_p)\stackrel{j^{\ast}}{\rightarrow}Ext_A^{3,p^nq+pq+1}(H^{\ast}M,H^{\ast}M)
\stackrel{i^{\ast}}{\rightarrow}\\
Ext_A^{3,p^nq+pq+1}(H^{\ast}M,\mathbb{Z}_p)=0
\end{array}$$
induced by (1.1). It is easy to get that the second group is zero.

At last consider the following exact sequence
$$\begin{array}{r}
0=Ext_A^{3,p^nq+(p+1)q+2}(H^{\ast}M,H^{\ast}M)\stackrel{i^{\prime}_{\ast}}{\rightarrow}
Ext_A^{3,p^nq+(p+1)q+2}(H^{\ast}K,H^{\ast}M)\\\stackrel{j^{\prime}_{\ast}}{\rightarrow}
Ext_A^{3,p^nq+pq+1}(H^{\ast}M,H^{\ast}M)=0
\end{array}$$ induced by (1.2). The
desired result follows.

{\protect\vspace{-15pt}\rightline{$\square$}}

{\bf Proposition 2.6} {\it Let $p\geq 5,n\geq3$, then
$$Ext_A^{2,p^nq+(p+1)q+1}(H^{\ast}K,\mathbb{Z}_p)
\cong \mathbb{Z}_p\{\beta_{\ast}i^{\prime}_{\ast}i_{\ast}(h_n)\},$$
where $\beta_{\ast}:Ext_A^{1,p^nq}(H^{\ast}K,\mathbb{Z}_p)\longrightarrow Ext_A^{2,p^nq+(p+1)q+1}(H^{\ast}K,\mathbb{Z}_p)$
is the connecting homomorphism induced by $\beta:\Sigma^{(p+1)q} K\longrightarrow
K$.}

{\it Proof.} Look at the exact sequence
$$\begin{array}{cl}
Ext_A^{1,p^nq+pq-1}(\mathbb{Z}_p,\mathbb{Z}_p)\stackrel{p_{\ast}}{\rightarrow}
Ext_A^{2,p^nq+pq}(\mathbb{Z}_p,\mathbb{Z}_p)&\stackrel{i_{\ast}}{\rightarrow}
Ext_A^{2,p^nq+pq}(H^{\ast}M,\mathbb{Z}_p)\\
&\stackrel{j_{\ast}}{\rightarrow}Ext_A^{2,p^nq+pq-1}(\mathbb{Z}_p,\mathbb{Z}_p)
\end{array}
$$
induced by (1.1). Since the first group and the fourth
group are zero, so the above $i_{\ast}$ is an
isomorphism. Thus we can see that $Ext_A^{2,p^nq+pq}(H^{\ast}M,\mathbb{Z}_p)\cong
\mathbb{Z}_p\{i_{\ast}(h_1h_n)\}$ by the
fact that $Ext_A^{2,p^nq+pq}(\mathbb{Z}_p,\mathbb{Z}_p)\cong
\mathbb{Z}_p\{h_1h_n\}$.

At last observe the following exact sequence
$$\begin{array}{r}
Ext_A^{2,p^nq+(p+1)q+1}(H^{\ast}M,\mathbb{Z}_p)\stackrel{i^{\prime}_{\ast}}{\rightarrow}
Ext_A^{2,p^nq+(p+1)q+1}(H^{\ast}K,\mathbb{Z}_p)\stackrel{j^{\prime}_{\ast}}{\rightarrow}\\
Ext_A^{2,p^nq+pq}(H^{\ast}M,\mathbb{Z}_p)\stackrel{\alpha_{\ast}}{\rightarrow}
Ext_A^{3,p^nq+(p+1)q+1}(H^{\ast}M,\mathbb{Z}_p)
\end{array}$$
induced by (1.2). Since the first group is zero by
$Ext_A^{2,p^nq+(p+1)q+r}(\mathbb{Z}_p,\mathbb{Z}_p)=0$ for $r=0,1$
(cf. [2]) and the fourth group is zero by
$Ext_A^{3,p^nq+(p+1)q+t}(\mathbb{Z}_p,\mathbb{Z}_p)=0$ for $t=0,1$
(cf. [6, Table 8.1]), then the above $j^{\prime}_{\ast}$ is an
isomorphism. Thus we can have that
$Ext_A^{2,p^nq+(p+1)q+1}(H^{\ast}K,\mathbb{Z}_p)$ has a unique
generator $\Delta$, which satisfies
$j^{\prime}_{\ast}(\Delta)=i_{\ast}(h_1h_n)$.  From [5, (5.4)], we
have that $j^{\prime}\beta i^{\prime}i\in [\sum^{pq-1}S,M]$ is
represented by $i_{\ast}(h_1)\in
Ext_A^{1,pq}(H^{\ast}M,\mathbb{Z}_p)$ in the Adams spectral
sequence. It follows that $(j^{\prime} \beta i^{\prime}
i)_{\ast}(h_n)=i_{\ast}(h_1h_n)=j^{\prime}_{\ast}(\Delta)$. Note the
fact that $j^{\prime}_{\ast}$ is an isomorphism. It is easy to get
that $\beta_{\ast}i^{\prime}_{\ast}i_{\ast}(h_n)=\Delta$. Therefore
this completes the proof of the proposition.

{\protect\vspace{-15pt}\rightline{$\square$}}

{\bf Proposition 2.7} {\it Let $p\geq 5,n\geq 3$, then
$$Ext_A^{2,p^nq+(p+1)q+1}(H^{\ast}K,H^{\ast}M)\cong \mathbb{Z}_p\{\beta_{\ast} i^{\prime}_{\ast}
(\tilde{h}_n)\}$$ where
$\tilde{h}_n\in Ext_A^{1,p^nq}(H^{\ast}M,H^{\ast}M)$ is the unique generator of
$Ext_A^{1,p^nq}(H^{\ast}M,H^{\ast}M)$ and satisfies
$i^{\ast}(\tilde{h}_n)=i_{\ast}(h_n)$.}

{\it Proof.} Consider the exact sequence
$$\begin{array}{cl}
Ext_A^{2,p^nq+pq+1}(\mathbb{Z}_p,\mathbb{Z}_p)\stackrel{i_{\ast}}{\rightarrow}Ext_A^{2,p^nq+pq+1}(H^{\ast}M,\mathbb{Z}_p)
&\stackrel{j_{\ast}}{\rightarrow}Ext_A^{2,p^nq+pq}(\mathbb{Z}_p,\mathbb{Z}_p)\\
&\stackrel{p_{\ast}}{\rightarrow}
Ext_A^{3,p^nq+pq+1}(\mathbb{Z}_p,\mathbb{Z}_p).
\end{array}$$
Since $Ext_A^{2,p^nq+pq+1}(\mathbb{Z}_p,\mathbb{Z}_p)=0$,
$imi_{\ast}=0$. Since
$Ext_A^{2,p^nq+pq}(\mathbb{Z}_p,\mathbb{Z}_p)=\mathbb{Z}_p\{h_1h_n\}$
and
$Ext_A^{3,p^nq+pq+1}(\mathbb{Z}_p,\mathbb{Z}_p)=\mathbb{Z}_p\{a_0h_1h_n\}$
(cf. [6, Table 8.1]), then the above $p_{\ast}$ is an isomorphism,
and then $imj_{\ast}=0$. Thus
$Ext_A^{2,p^nq+pq+1}(H^{\ast}M,\mathbb{Z}_p)=0$.

Look at the exact sequence
$$\begin{array}{r}Ext_A^{2,p^nq+(p+1)q+2}(H^{\ast}M,\mathbb{Z}_p)\stackrel{i^{\prime}_{\ast}}{\rightarrow}
Ext_A^{2,p^nq+(p+1)q+2}(H^{\ast}K,\mathbb{Z}_p)\stackrel{j^{\prime}_{\ast}}{\rightarrow}\\
Ext_A^{2,p^nq+pq+1}(H^{\ast}M,\mathbb{Z}_p)=0
\end{array}$$ induced by (1.2). Since we know that
$Ext_A^{2,p^nq+(p+1)q+2}(H^{\ast}M.\mathbb{Z}_p)=0$ by the
facts that $Ext_A^{2,p^nq+(p+1)q+r}(\mathbb{Z}_p,\mathbb{Z}_p)=0$
for $r=1,2$ and $Ext_A^{2,p^nq+pq+1}(H^{\ast}M,\mathbb{Z}_p)=0$, then
$Ext_A^{2,p^nq+(p+1)q+2}(H^{\ast}K,\mathbb{Z}_p)=0$.

Observe the following exact sequence
$$Ext_A^{3,p^nq+(p+1)q+2}(H^{\ast}M,\mathbb{Z}_p)\stackrel{i^{\prime}_{\ast}}{\rightarrow}
Ext_A^{3,p^nq+(p+1)q+2}(H^{\ast}K,\mathbb{Z}_p)\stackrel{j^{\prime}_{\ast}}{\rightarrow}
Ext_A^{3,p^nq+pq+1}(H^{\ast}M,\mathbb{Z}_p)$$ induced by (1.2). From the
proof of Proposition 2.5 we know that the first group and the
third group are zero. So the middle group is zero.

At last, look at the following exact sequence
$$\begin{array}{r}
Ext_A^{2,p^nq+(p+1)q+2}(H^{\ast}K,\mathbb{Z}_p)\stackrel{j^{\ast}}{\rightarrow}
Ext_A^{2,p^nq+(p+1)q+1}(H^{\ast}K,H^{\ast}M)\stackrel{i^{\ast}}{\rightarrow}\\
Ext_A^{2,p^nq+(p+1)q+1}(H^{\ast}K,\mathbb{Z}_p)
\stackrel{p^{\ast}}{\rightarrow}Ext_A^{3,p^nq+(p+1)q+2}(H^{\ast}K,\mathbb{Z}_p)
\end{array}$$
induced by (1.1). Since the first group and the fourth group are
zero, then the above $i^{\ast}$ is an isomorphism. Notice that
$Ext_A^{2,p^nq+(p+1)q+1}(H^{\ast}K,\mathbb{Z}_p)\cong
\mathbb{Z}_p\{\beta_{\ast}i^{\prime}_{\ast}i_{\ast}(h_n)\}$. Thus we
can easily have that there exists an element $\bar{\bar{\Delta}}\in
Ext_A^{2,p^nq+(p+1)q+1}(H^{\ast}K,H^{\ast}M)$ such that
$Ext_A^{2,p^nq+(p+1)q+1}(H^{\ast}K,H^{\ast}M)\cong
\mathbb{Z}_p\{\bar{\bar{\Delta}}\}$ and
$i^{\ast}(\bar{\bar{\Delta}})
=\beta_{\ast}i^{\prime}_{\ast}i_{\ast}(h_n)$. Since
$i^{\ast}(\tilde{h}_n)=i_{\ast}(h_n)$, $i^{\ast}(\bar{\bar{\Delta}})
=\beta_{\ast}i^{\prime}_{\ast}i_{\ast}(h_n)=\beta_{\ast}i^{\prime}_{\ast}
i^{\ast}(\tilde{h}_n)=i^{\ast}\beta_{\ast}i^{\prime}_{\ast}(\tilde{h}_n)$.
Thus we have that
$\bar{\bar{\Delta}}=\beta_{\ast}i^{\prime}_{\ast}(\tilde{h}_n)$ by
the fact that $i^{\ast}$ is an isomorphism. Thus
this completes the proof of the proposition.\\

{\protect\vspace{-15pt}\rightline{$\square$}}

\noindent{\bf\large 3 Proofs of the main theorems }
\\
\\
\indent
Let
$$
\begin{array}{cccccc}
\cdots\stackrel{\bar{a}_2}{\rightarrow}&\Sigma^{-2}E_2&
\stackrel{\bar{a}_1}{\rightarrow}&\Sigma^{-1}E_1&\stackrel{\bar{a}_0}
{\rightarrow}&S\\
&\downarrow \bar{b}_2& &\downarrow \bar{b}_1& &\downarrow \bar{b}_0\\
& \Sigma^{-2}KG_2& &\Sigma^{-1}KG_1& &KG_0=K\mathbb{Z}_p
\end{array}
\eqno{(3.1)}
$$
be the minimal Adams resolution of $S$ satisfying the following.

(1) $E_s\stackrel{\bar{b}_s}{\rightarrow}KG_s\stackrel{\bar{c}_s}{\rightarrow}E_{s+1}
\stackrel{\bar{a}_s}{\rightarrow}\Sigma E_s$ are cofibrations
for all $s\geq 0$ which induce short exact sequences in
$\mathbb{Z}_p$-cohomology.

(2) $KG_s$ is a wedge sum of suspensions of Eilenberg-Maclane
spectra of type $K\mathbb{Z}_p$.

(3) $\pi_t KG_s$ are the $E_1^{s,t}$-terms, $(\bar{b}_s
\bar{c}_{s-1})_{\ast}: \pi_t KG_{s-1}\rightarrow \pi_t KG_s$ are the
$d^{s-1,t}_1$-differentials of the Adams spectral sequence and
$\pi_t KG_s\cong Ext_A^{s,t}(\mathbb{Z}_p,\mathbb{Z}_p)$ (cf. [9,
p.180]). Then
$$
\begin{array}{cccccl}
\cdots\stackrel{\bar{a}_2\wedge 1_W}{\rightarrow}&\Sigma^{-2}E_2\wedge W&
\stackrel{\bar{a}_1 \wedge 1_W}{\rightarrow}&\Sigma^{-1}E_1 \wedge W&
\stackrel{\bar{a}_0 \wedge 1_W}
{\rightarrow}&W\\
&\downarrow \bar{b}_2 \wedge 1_W& &\downarrow \bar{b}_1 \wedge 1_W&
&\downarrow \bar{b}_0 \wedge 1_W\\
& \Sigma^{-2}KG_2 \wedge W& &\Sigma^{-1}KG_1 \wedge W& &KG_0\wedge
W
\end{array}\eqno{(3.2)}
$$
is an Adams resolution of arbitrary finite spectrum W.

From [10, p.204-206], the Moore spectrum $M$ is a commutative ring
spectrum with multiplication $m_M:M\wedge M\rightarrow M$ and there
is $\bar{m}_M:\Sigma M\rightarrow M\wedge M$ such that
 $$\begin{array}{ll}
 m_M(i\wedge 1_M)=1_M,&(j\wedge 1_M)\bar{m}_M=1_M,\\
m_M\bar{m}_M=0,&(i\wedge 1_M)m_M+\bar{m}_M(j\wedge 1_M)=1_{M\wedge
M},\\
m_M T=-m_M,&T\bar{m}_M=\bar{m}_M,\\
m_M(1_M\wedge i)=-1_M,&(1_M\wedge j)\bar{m}_M=1_M,
\end{array}
$$
where $T: M\wedge M\longrightarrow M\wedge M$ is the switching map.

A spectrum $X$ is called an $M$-module spectrum if $p\wedge
1_X=0$, and consequently, the cofibration
$X\stackrel{p\wedge 1_X}{\longrightarrow}X\stackrel{i\wedge 1_X}
{\longrightarrow}M\wedge X\stackrel{j\wedge 1_X}{\longrightarrow}\Sigma X$
split, i.e., there is a homotopy equivalence
$M\wedge X=X\bigvee \Sigma X$ and there are maps $m_X: M\wedge X\longrightarrow
X$, $\bar{m}_X: \Sigma X\longrightarrow M\wedge X$ satisfying
$m_X(i\wedge 1_X)=1_X$, $(j\wedge 1_X)\bar{m}_X=1_X$,
$m_X\bar{m}_X=0$ and $\bar{m}_X(j\wedge 1_X)+(i\wedge 1_X)m_X=1_{M\wedge X}$.
The $M$-module actions $m_X$, $\bar{m}_X$ are called associative if
$m_X(1_M\wedge m_X)=-m_X(m_X\wedge 1_X)$
and $(1_M\wedge\bar{m}_M)\bar{m}_X=
(\bar{m}_M\wedge 1_X)\bar{m}_X$.

Let $X$ and $X^{\prime}$ be $M$-module spectra. Then we define a
homomorphism $d:[\Sigma^s X^{\prime}, X]\rightarrow
[\Sigma^{s+1}X^{\prime}, X]$ by $d(f)=m_X(1_M\wedge
f)\bar{m}_{X^{\prime}}$ for $f\in [\Sigma^s X^{\prime}, X]$. This
operation $d$ is called a derivation (of maps between $M$-module
spectra) which has the following properties:

{\bf Lemma} ${{\bf 3.1}}^{\hbox{[10, Theorem 2.2]}}$(1) {\it $d$ is
a derivative: $d(fg)=fd(g)+(-1)^{|g|}d(f)g$ for $f\in [\Sigma^s
X^{\prime}, X]$, $g\in [\Sigma^{t}X^{\prime\prime}, X^{\prime}]$,
where $X$, $X^{\prime}$, $X^{\prime\prime}$ are $M$-module spectra.

(2) Let $W^{\prime}$, $W$ be arbitrary spectra and $h\in [\Sigma^{r}W^{\prime},
W]$. Then $$d(h\wedge f)=(-1)^{|h|}h\wedge d(f)$$ for $f\in [\Sigma^{s}X^{\prime},
X]$.

(3) $d^2=0:[\Sigma^{s}X^{\prime}, X]\longrightarrow [\Sigma^{s+2}X^{\prime},
X]$ for associative spectra $X^{\prime}$, $X$.}

{\protect\vspace{-15pt}\rightline{$\square$}}

 From [10,
(3.4)], $K$ is an $M$-module spectrum, i.e., there are $M$-module
actions $m_K: K\wedge M\longrightarrow K$, $\bar{m}_K:\Sigma
K\longrightarrow K\wedge M$ satisfying $$m_K(1_K\wedge i)=1_K,
(1_K\wedge j)\bar{m}_K=1_K, m_K\bar{m}_K=0, (1_K\wedge
i)m_K+(1_K\wedge j)\bar{m}_K=1_{K\wedge M}.$$ Moreover, from [10,
(2.6)] and [10, (3.7)] we have that $d(ij)=-1_M$, $d(\alpha)=0$,
$d(i^{\prime})=0$, $d(j^{\prime})=0$ and $d(\beta)=0$.

{\bf Remark 3.1}\quad {\it In this paper, all the notations are the
same as those of [7].}

Let $L$ be the cofiber of $\alpha_1=j\alpha
i:\Sigma^{q-1}S\rightarrow S$ and $K^{\prime}$ be the cofiber of
$\alpha i:\Sigma^{q}S\rightarrow M$ given by the following two
cofibrations:
$$
\Sigma^{q-1}S\stackrel{\alpha_1}{\rightarrow}S\stackrel
{i^{\prime\prime}}{\rightarrow}L\stackrel{j^{\prime\prime}}{\rightarrow}\Sigma^{q}S
(see [7, (2.3)]), \eqno{(3.3)}
$$
$$
\Sigma^{q}S\stackrel{\alpha
i}{\rightarrow}M\stackrel{\upsilon}{\rightarrow}K^{\prime}\stackrel{y}
{\rightarrow}\Sigma^{q+1}S (see [7, (2.4)]. \eqno{(3.4)}
$$
Let $\alpha^{\prime}=\alpha_1 \wedge 1_K$. Consider the following
two commutative diagram of $3\times 3$ in the stable homotopy
category:

$${\small \begin{array}{clclccclc} \Sigma
M&&\stackrel{\upsilon}{\longrightarrow}&&\Sigma
K^{\prime}&&\stackrel{1_{K^{\prime}}\wedge
p}{\longrightarrow}&&\Sigma K^{\prime}\\
&\searrow (\upsilon\wedge 1_M)\bar{m}_M& &\nearrow
1_{K^{\prime}}\wedge j& &\searrow y& & \nearrow z&\\
& &K^{\prime}\wedge M& & & &\Sigma^{q+2}S& & \\
&\nearrow 1_{K^{\prime}}\wedge i& &\searrow\pi& & \nearrow
jj^{\prime}&
&\searrow\alpha i&\\
K^{\prime}& &\stackrel{x}{\longrightarrow}& &K&
&\stackrel{j^{\prime}\alpha^{\prime}}{\longrightarrow}&& \Sigma^2 M
\end{array}
}$$ and
$${\small \begin{array}{clclccclc}
M&&\stackrel{(i^{\prime\prime}\wedge
1_K)i^{\prime}}{\longrightarrow}&&L\wedge
K&&\stackrel{j^{\prime\prime}\wedge
1_K}{\longrightarrow}&&\Sigma^q K\\
&\searrow i^{\prime}& &\nearrow
i^{\prime\prime}\wedge 1_K& &\searrow \bar{r}& & \nearrow \pi&\\
& &K& & & &\Sigma^{q}K^{\prime}\wedge M& & \\
&\nearrow \alpha^{\prime}& &\searrow j^{\prime}& & \nearrow
(\upsilon\wedge 1_M)\bar{m}_M&
&\searrow \varepsilon&\\
\Sigma^{q-1}K&
&\stackrel{j^{\prime}\alpha^{\prime}}{\longrightarrow}&
&\Sigma^{q+1}M& &\stackrel{\alpha}{\longrightarrow}&& \Sigma M
\end{array}}
$$

By the above two commutative diagram of $3\times 3$ in the stable
homotopy category, we easily have the following two lemmas.

{\bf Lemma 3.2}\quad {\it There exist three cofibrations
$$
 K^{\prime}\stackrel{x}{\longrightarrow}K\stackrel{jj^{\prime}}{\longrightarrow}
 \Sigma^{q+2} S  \stackrel{z}{\longrightarrow}\Sigma K^{\prime}
,\eqno{(3.5)}
$$
$$
 \Sigma^{-1}K\stackrel{j^{\prime}\alpha^{\prime}}{\longrightarrow}\Sigma
 M\stackrel{(\upsilon\wedge
 1_{M}){\bar{m}}_M}{\longrightarrow}K^{\prime}\wedge
 M\stackrel{\pi}{\longrightarrow}K
,\eqno{(3.6)}
$$
$$
M\stackrel{(i^{\prime\prime}\wedge
1_K)i^{\prime}}{\longrightarrow}L\wedge
K\stackrel{\bar{r}}{\longrightarrow}\Sigma^{q}K^{\prime}\wedge
M\stackrel{\varepsilon}{\longrightarrow}\Sigma M.
\eqno{(3.7)}
$$}

{\protect\vspace{-15pt}\rightline{$\square$}}

 {\bf Lemma 3.3}\quad
{\it $\varepsilon(\upsilon\wedge 1_M){\bar{m}_M} =\alpha$,
$\bar{r}(i^{\prime\prime}\wedge 1_K)= (\upsilon\wedge
1_M)\bar{m}_Mj^{\prime}$, $\pi\bar{r} =j^{\prime\prime}\wedge 1_K$,
$\varepsilon(1_{K^{\prime}}\wedge i)vj^{\prime}=
-2j^{\prime}\alpha^{\prime}.$}

{\protect\vspace{-15pt}\rightline{$\square$}}

From [11, p.434], there are $\bar{\triangle}\in
[\Sigma^{-1-1}L\wedge K,K]$ and $\tilde{\Delta}\in
[\Sigma^{-1}K,L\wedge K]$ satisfying
$\bar{\Delta}(i^{\prime\prime}\wedge 1_{K})=(j^{\prime\prime}\wedge
1_{K})\tilde{\Delta}=i^{\prime}j^{\prime}\in [\Sigma^{-q-1}K,K]$ and
$jj^{\prime}\bar{\Delta}=0$. From [6, p.484], there is
$\bar{\Delta}_{K^{\prime}}\in [\Sigma^{-q-1}L\wedge K,K^{\prime}]$
such that $\bar{\Delta}_{K^{\prime}}(i^{\prime\prime}\wedge
1_{K})=\upsilon j^{\prime}\in [\Sigma^{-q-1} K,K^{\prime}]$ and
$\bar{\Delta}(i^{\prime\prime}\wedge 1_K)=(j^{\prime\prime}\wedge
1_{K})\tilde{\Delta}=i^{\prime}j^{\prime}$.

{\bf Lemma 3.4}\quad {\it
$\bar{\Delta}_{K^{\prime}}=(1_{K^{\prime}}\wedge j)\bar{r}$.}

{\it Proof.} From Lemma 3.3 we have
$$(1_{K^{\prime}}\wedge j)\bar{r}(i^{\prime\prime}\wedge 1_K)=(1_K^{\prime}\wedge j)(\upsilon
\wedge 1_M)\bar{m}_{M}j^{\prime}=(\upsilon\wedge 1_{S^0})(1_M\wedge j)\bar{m}_M j^{\prime}
=\upsilon j^{\prime}=\bar{\Delta}_{K^{\prime}}(i^{\prime\prime}\wedge 1_K),$$
which shows that
$(1_{K^{\prime}}\wedge j)\bar{r}=\bar{\Delta}_{K^{\prime}}+\bar{g}(j^{\prime\prime}\wedge 1_K)$
for some $\bar{g}\in [K,\Sigma K^{\prime}]$

Consider the exact sequence induced by
(3.4)
$$[K,\Sigma^{q+1}S]\stackrel{(\alpha i)_{\ast}}{\rightarrow}[K,\Sigma M]\stackrel{{\upsilon}_{\ast}}
{\rightarrow}[K,\Sigma
K^{\prime}]\stackrel{y_{\ast}}{\rightarrow}[K,
\Sigma^{q+2}S]\stackrel{(\alpha
i)_{\ast}}{\rightarrow}[K,\Sigma^{2}M].$$ From the proof of [7,
Proposition 2.18], we know that $[K,\Sigma M]=0$. So $im
\upsilon_{\ast}=0$. Since $[K,\Sigma^{q+2}S]\cong
\mathbb{Z}_p\{jj^{\prime}\}$, $(\alpha i)_{\ast}(jj^{\prime})=\alpha
ijj^{\prime}\not=0$. Thus we have $im y_{\ast}=0$ and $[K,\Sigma
K^{\prime}]=0$. Then we have $(1_{K^{\prime}}\wedge
j)\bar{r}=\bar{\Delta}_{K^{\prime}}$.

{\protect\vspace{-15pt}\rightline{$\square$}}

{\bf Lemma 3.5}[7, lemma 3.3 and (3.4)] {\it Let $p\geq 5, n\geq 3$,
then there exists an element $\eta^{\prime}_{n,2}\in [\Sigma^{p^nq+q}K,
E_2\wedge K]$ such that $$(\bar{b}_2\wedge
1_K)\eta^{\prime}_{n,2}=h_0h_n\wedge 1_K\in [\Sigma^{p^nq+q}K,
KG_2\wedge K], (1_{E_2}\wedge
\alpha^{\prime})\eta^{\prime}_{n,2}=0,$$ where
$h_0h_n\in\pi_{p^nq+q}KG_2\cong
Ext_A^{2,p^nq+q}(\mathbb{Z}_p,\mathbb{Z}_p)$ and
$\alpha^{\prime}=j\alpha i\wedge 1_K\in [\Sigma^{q-1}K, K]$. There
also exists an element $f_2\in [\Sigma^{p^nq+(p+2)q+3}M, E_5\wedge L\wedge K]$
such that $$(1_{E_2}\wedge (i^{\prime\prime}\wedge
1_K)\beta){\eta^{\prime}}_{n,2}i^{\prime}
=({\bar{a}}_2{\bar{a}}_3{\bar{a}}_4\wedge 1_{L\wedge K})f_2.$$}

{\protect\vspace{-15pt}\rightline{$\square$}}

 {\bf Corollary
3.1}\quad {\it For $f_2\in [\Sigma^{p^nq+(p+2)q+3}M, E_5\wedge
L\wedge K]$ which is given in Lemma 3.5, we have
$$
(1_{E_4}\wedge \varepsilon(1_{K^{\prime}}\wedge
ij))({\bar{a}}_4\wedge 1_{K^{\prime}\wedge M})(1_{E_5}\wedge
\bar{r})d(f_2 ij)=0.\eqno{(3.8)}
$$}
\indent{\it Proof.} From [7], we have [7, (3.6)] that
$(\bar{a}_4\wedge 1_M)(1_{E_5}\wedge
\varepsilon(1_{K^{\prime}}\wedge
i){\bar{\Delta}}_{K^{\prime}})d(f_2ij)=0.$ Here $f_2$ is given in
[7, (3.4)].

By [10, (1.7)], we have that $(\bar{a}_4\wedge 1_M)(1_{E_5}\wedge
\varepsilon(1_{K^{\prime}}\wedge i))(1_{E_5}\wedge
{\bar{\Delta}}_{K^{\prime}})d(f_2ij)=0.$ By lemma 3.4, we have
$$(\bar{a}_4\wedge 1_M)(1_{E_5}\wedge
\varepsilon(1_{K^{\prime}}\wedge i))(1_{E_5}\wedge
(1_{K^{\prime}}\wedge j)\bar{r})d(f_2ij)=0.$$ By [10, (1.7)], it
follows that $$(\bar{a}_4\wedge 1_M)(1_{E_5}\wedge
\varepsilon(1_{K^{\prime}}\wedge i))(1_{E_5}\wedge
(1_{K^{\prime}}\wedge j))(1_{E_5}\wedge\bar{r})d(f_2ij)=0.$$ Thus
$$(\bar{a}_4\wedge 1_M)(1_{E_5}\wedge
\varepsilon(1_{K^{\prime}}\wedge
ij))(1_{E_5}\wedge\bar{r})d(f_2ij)=0.$$ By [10, (1.7)], the
corollary follows.

{\protect\vspace{-15pt}\rightline{$\square$}}

Let $W$ be the cofibre of $\varepsilon(1_{K^{\prime}}\wedge
ij):\Sigma^{q-2}K^{\prime}\wedge M\longrightarrow M$ given by the
cofibration
$$
\Sigma^{q-2}K^{\prime}\wedge M\stackrel{\varepsilon(1_{K^{\prime}}\wedge ij)}{\longrightarrow}
M\stackrel{w_4}{\longrightarrow}
W\stackrel{u_4}{\longrightarrow}\Sigma^{q-1}K^{\prime}\wedge M.\eqno{(3.9)}
$$

{\bf Lemma 3.6} {\it There exists $f^{\prime}\in
[\Sigma^{p^nq+(p+2)q+1}M,E_4\wedge W]$ such that $
({\bar{a}}_2{\bar{a}}_3\wedge 1_K)(1_{E_4}\wedge \pi
u_4)f^{\prime}=0. $}

{\it Proof.}\quad By (3.8) and (3.9), we have that
$$
({\bar{a}}_4\wedge 1_{K^{\prime}\wedge M})(1_{E_5}\wedge
\bar{r})d(f_2ij)=(1_{E_4}\wedge u_4)f^{\prime}
\eqno{(3.10)}
$$
with $f^{\prime}\in [\Sigma^{p^nq+(p+2)q+1}M,E_4\wedge W]$ and by
composing $({\bar{a}}_2{\bar{a}}_3\wedge 1_K)(1_{E_4}\wedge\pi)$
on (3.10) we have
$$
({\bar{a}}_2{\bar{a}}_3\wedge 1_K)(1_{E_4}\wedge \pi
u_4)f^{\prime}=({\bar{a}}_2{\bar{a}}_3{\bar{a}}_4\wedge 1_K)(1_{E_5}\wedge \pi
\bar{r})d(f_2ij).\eqno{(3.11)}
$$

By composing $ij$ on [7, (3.4)], we have
$$
(1_{E_2}\wedge (i^{\prime\prime}\wedge
1_K)\beta){\eta^{\prime}}_{n,2}i^{\prime}ij
=({\bar{a}}_2{\bar{a}}_3{\bar{a}}_4\wedge 1_{L\wedge
K})f_2ij
\eqno{(3.12)}
$$
with ${\eta^{\prime}}_{n,2}\in [\Sigma^{p^nq+q}K, E_2\wedge K].$

Notice that $d(1_K)=0$ and $d(\beta)=0$. Then by applying the
derivation $d$ on (3.12) we have
$$
(1_{E_2}\wedge (i^{\prime\prime}\wedge
1_K)\beta)d({\eta^{\prime}}_{n,2}i^{\prime}ij)
=({\bar{a}}_2{\bar{a}}_3{\bar{a}}_4\wedge 1_{L\wedge
K})d(f_2ij)
.\eqno{(3.13)}
$$

Notice that $\pi\bar{r}=j^{\prime\prime}\wedge 1_K$. By composing $(1_{E_2}\wedge
\pi\bar{r})$ on (3.13) we have
$$
({\bar{a}}_2{\bar{a}}_3{\bar{a}}_4\wedge 1_K)(1_{E_5}\wedge
\pi\bar{r})d(f_2ij)=0
\eqno{(3.14)}
$$
and by (3.11), (3.14) we get
$$
({\bar{a}}_2{\bar{a}}_3\wedge 1_K)(1_{E_4}\wedge \pi
u_4)f^{\prime}=0.
\eqno{(3.15)}
$$
Thus the Lemma is proved.

{\protect\vspace{-15pt}\rightline{$\square$}}

Let $U$ be the cofibre of $\pi u_4:W\longrightarrow \Sigma^{q-1}K$
given by the cofibration
$$
W\stackrel{\pi
u_4}{\longrightarrow}\Sigma^{q-1}K\stackrel{w_5}{\longrightarrow}U\stackrel{u_5}{\longrightarrow}\Sigma
W.
\eqno{(3.16)
}
$$

{\bf Lemma 3.7}\quad {\it $w_5$ induces zero homomorphism in
$\mathbb{Z}_p$-cohomology.}

{\it Proof.} Consider the following homomorphism induced by $w_5$:
$$w_5^{\ast}: H^{\ast}U\longrightarrow H^{\ast+q-1}K.$$
From the celluar structures of $U$ and $K$, we can have that
$$H^{t}K=\left\{\begin{array}{ll} \mathbb{Z}_p,& t=0,1,q+1,q+2;\\
0,&others,
\end{array}\right.$$ and the top cell of $U$ has degree $2q+1$.
It easily follows that $w_5^{\ast}$ must be a zero homomorphism in
$\mathbb{Z}_p$-cohomology.

{\protect\vspace{-15pt}\rightline{$\square$}}

{\bf Lemma 3.8}\quad {\it There exist three homotopy elements $f^{\prime}_2\in
[\Sigma^{p^nq+(p+2)q}M,E_2\wedge U]$, $f_3^{\prime}\in
[\Sigma^{p^nq+(p+2)q+1}M, E_3\wedge U]$ and $g_2\in
[\Sigma^{p^nq+(p+2)q}M, KG_2\wedge W]$ such that
$$({\bar{a}}_2{\bar{a}}_3\wedge 1_W)f^{\prime}=(1_{E_2}\wedge
u_5)f_2^{\prime}, f_2^{\prime}=({\bar{a}}_2\wedge
1_U)f_3^{\prime}$$ and $$(1_{E_3}\wedge u_4)({\bar{a}}_3\wedge
1_W)f^{\prime}=-(1_{E_3}\wedge u_4 u_5)f_3^{\prime}+(1_{E_3}\wedge
u_4)({\bar{c}}_2\wedge 1_W)g_2.$$}

{\it Proof.}\quad From (3.15) and (3.16), we have that
$$
({\bar{a}}_2{\bar{a}}_3\wedge 1_W)f^{\prime}=(1_{E_2}\wedge
u_5)f_2^{\prime}
\eqno{(3.17)}
$$
with $f^{\prime}_2\in [\Sigma^{p^nq+(p+2)q}M,E_2\wedge U]$.

By
(3.17) and (3.2) we have $({\bar{b}}_2\wedge 1_W)(1_{E_2}\wedge
u_5)f_2^{\prime}=({\bar{b}}_2\wedge 1_W)
({\bar{a}}_2{\bar{a}}_3\wedge 1_W)f^{\prime}=0$. Thus
$$
(1_{KG_2}\wedge u_5)({\bar{b}}_2\wedge
1_U)f_2^{\prime}=0.\eqno{(3.18)}
$$
By (3.18) ,(3.16) and the fact that $w_5$ induces zero homomorphism
in $\mathbb{Z}_p$-cohomology(see lemma 3.7), we have
$$
({\bar{b}}_2\wedge 1_U)f_2^{\prime}=(1_{KG_2}\wedge
w_5)g=0
\eqno{(3.19)}
$$
with $g\in [\Sigma^{p^nq+(p+1)q+1}M,KG_2\wedge K]$, so by (3.2) we obtain
$$
f_2^{\prime}=({\bar{a}}_2\wedge 1_U)f_3^{\prime}
\eqno{(3.20)}
$$
with $f_3^{\prime}\in [\Sigma^{p^nq+(p+2)q+1}M, E_3\wedge U]$. By
[10, (1.7)], from (3.20) and (3.17) we have
$({\bar{a}}_2{\bar{a}}_3\wedge 1_W)f^{\prime}= -({\bar{a}}_2\wedge
1_W)(1_{E_3}\wedge u_5)f_3^{\prime}$. Then we have
$$
({\bar{a}}_3\wedge 1_W)f^{\prime}=-(1_{E_3}\wedge
u_5)f_3^{\prime}+({\bar{c}}_2\wedge 1_W)g_2
\eqno{(3.21)}
$$
with $g_2\in [\Sigma^{p^nq+(p+2)q}M, KG_2\wedge W]$. By composing
$(1_{E_3}\wedge u_4)$ on (3.21), we have
$$
(1_{E_3}\wedge u_4)({\bar{a}}_3\wedge 1_W)f^{\prime}=-(1_{E_3}\wedge u_4
u_5)f_3^{\prime}+(1_{E_3}\wedge u_4)({\bar{c}}_2\wedge 1_W)g_2.
\eqno{(3.22)}
$$
We finishes the proof of the lemma.

{\protect\vspace{-15pt}\rightline{$\square$}}

{\bf Lemma 3.9}\quad {\it the cofibre of
$\varepsilon(1_{K^{\prime}}\wedge i)\upsilon: \Sigma^q
M\longrightarrow \Sigma M$ is $U$ given by the cofibration
$$
\Sigma^q M\stackrel{\varepsilon(1_{K^{\prime}}\wedge
i)\upsilon}{\longrightarrow}\Sigma M
\stackrel{w_6}{\longrightarrow}U\stackrel{u_6}{\longrightarrow}\Sigma^{q+1}M.
\eqno{(3.23)}
$$
There exist two relations that $$u_4 u_5=(\upsilon\wedge
1_M){\bar{m}}_M u_6$$ and $$\varepsilon(1_{K^{\prime}}\wedge
ij)(\upsilon\wedge 1_M)\bar{m}_M =\varepsilon(1_{K^{\prime}}\wedge
i)\upsilon.$$}

{\it Proof.}\quad By the three cofibrations (3.6), (3.9), and
(3.16), we can get the following commutative diagram (3.24) of
$3\times 3$ lemma in stable homotopy category(cf. [12, p. 292-293]).
$${\tiny
\begin{array}{rllcccclll}
W &\stackrel{\pi u_4}{\longrightarrow}& &\Sigma^{q-1}K&
&\stackrel{j^{\prime}\alpha^{\prime}}{\longrightarrow}&
&\Sigma^{q+1}M& & \\
u_4\searrow& &\pi\nearrow& &\searrow w_5&
&\nearrow u_6& &\searrow\upsilon\wedge 1_M\bar{m}_M& \\
  &\Sigma^{q-1}K^{\prime}\wedge M& & & &U& & &
 &\Sigma^{q}K^{\prime}\wedge M.\\
 \upsilon\wedge 1_M\bar{m}_M\nearrow& &\searrow\varepsilon(1_{K^{\prime}}\wedge
  ij)& &\nearrow w_6& &\searrow u_5& &\nearrow u_4&\\
  \Sigma^{q}M &\stackrel{\varepsilon(1_{K^{\prime}}\wedge
  i)\upsilon}{\longrightarrow}& &\Sigma M&
  &\stackrel{w_4}{\longrightarrow}& &\Sigma W& &
  \end{array}}$$
By the commutative diagram (3.24), Lemma 3.8 follows.

{\protect\vspace{-15pt}\rightline{$\square$}}

{\bf Lemma 3.10}\quad {\it With notation as above. We have
$$\begin{array}{cl} &({\bar{a}}_3{\bar{a}}_4\wedge
1_{K^{\prime}\wedge M})(1_{E_5}\wedge\bar{r})d(f_2ij)
\\=&(1_{E_3}\wedge (\upsilon\wedge 1_M){\bar{m}}_M u_6)f_3^{\prime}-({\bar{c}}_2
\wedge 1_{K^{\prime}\wedge M})
 (1_{KG_2}\wedge u_4)g_2.\end{array}.$$}

 {\it Proof.}\quad By
(3.22), [10, (1.7)] and the relation $u_4 u_5=(\upsilon\wedge
1_M){\bar{m}}_M u_6$ (see Lemma 3.9), we have
$$
({\bar{a}}_3\wedge 1_{K^{\prime}\wedge M})(1_{E_4}\wedge
u_4)f^{\prime}
=(1_{E_3}\wedge (\upsilon\wedge 1_M){\bar{m}}_M u_6)f_3^{\prime}-({\bar{c}}_2\wedge 1_{K^{\prime}\wedge M})
(1_{KG_2}\wedge u_4)g_2.
\eqno{(3.25)}
$$
By composing $({\bar{a}}_3\wedge 1_{K^{\prime}\wedge M})$ on (3.10),
we have
$$
({\bar{a}}_3\wedge 1_{K^{\prime}\wedge M})(1_{E_4}\wedge
u_4)f^{\prime}=
({\bar{a}}_3{\bar{a}}_4\wedge 1_{K^{\prime}\wedge M})(1_{E_5}\wedge
\bar{r})d(f_2ij).\eqno{(3.26)}
$$
Combining (3.25) and (3.26) yields
$$
\begin{array}{cl}
&({\bar{a}}_3{\bar{a}}_4\wedge 1_{K^{\prime}\wedge M})(1_{E_5}\wedge\bar{r})d(f_2ij)
\\=&(1_{E_3}\wedge (\upsilon\wedge 1_M){\bar{m}}_M u_6)f_3^{\prime}-({\bar{c}}_2
\wedge 1_{K^{\prime}\wedge M})
 (1_{KG_2}\wedge u_4)g_2.\end{array}
 \eqno{(3.27)}
 $$
 Thus we complete the proof of this lemma.

 {\protect\vspace{-15pt}\rightline{$\square$}}

{\bf Lemma 3.11}\quad {\it There exist two elements $f_4^{\prime}\in
[\Sigma^{p^nq+(p+2)q+1}M,E_3\wedge K]$ and $f_5^{\prime}\in
[\Sigma^{p^nq+(p+2)q+1}M, E_3\wedge K^{\prime}\wedge M]$ such that
$$ (1_{E_3}\wedge u_6)f_3^{\prime}=(1_{E_3}\wedge
j^{\prime})f_4^{\prime}$$ and $$
f_4^{\prime}=(1_{E_3}\wedge\pi)f_5^{\prime}.$$}

{\it Proof.}\quad By Lemma 3.3, we have
$\alpha=\varepsilon(\upsilon\wedge 1_M){\bar{m}}_M$. Then
$$
\begin{array}{ll}
  &(1_{E_3}\wedge\alpha u_6)f_3^{\prime} \\
= &(1_{E_3}\wedge \varepsilon(\upsilon\wedge
1_M)\bar{m}_M u_6)f_3^{\prime},\quad\hbox{since $u_4 u_5
 =(\upsilon\wedge 1_M)\bar{m}_M u_6$}\\
 = & (1_{E_3}\wedge \varepsilon u_4 u_5)f_3^{\prime}\\
= &(1_{E_3}\wedge \varepsilon)(1_{E_3}\wedge u_4 u_5)f_3^{\prime},\quad\hbox{by (3.22)}
\\= & (1_{E_3}\wedge \varepsilon)[(1_{E_3}\wedge u_4)({\bar{c}}_2\wedge 1_W)g_2-(1_{E_3}
\wedge u_4)({\bar{a}}_3\wedge 1_W)f^{\prime}]\\
= &(1_{E_3}\wedge\varepsilon u_4)(\bar{c}_2\wedge 1_W)g_2-(1_{E_3}
\wedge \varepsilon u_4)(\bar{a}_3\wedge 1_W)f^{\prime} \\
= & (\bar{c}_2\wedge 1_M)(1_{KG_2}
\wedge \varepsilon u_4)g_2-(1_{E_3}\wedge \varepsilon u_4)
(\bar{a}_3\wedge 1_W)f^{\prime},\quad\hbox{since $1_{KG_2}\wedge\varepsilon\simeq 0$} \\
= &-(1_{E_3}\wedge \varepsilon u_4)(\bar{a}_3
\wedge 1_W)f^{\prime}\\
 = &(1_{E_3}\wedge \varepsilon)(\bar{a}_3\wedge 1_{K^{\prime}\wedge M})
 (1_{E_4}\wedge u_4)f^{\prime},\quad\hbox{by (3.26)} \\
  = & (1_{E_3}\wedge \varepsilon)(\bar{a}_3\wedge 1_{K^{\prime}\wedge M})
  (\bar{a}_4\wedge 1_{K^{\prime}\wedge M})(1_{E_5}\wedge\bar{r})d(f_2ij)\\
= & (\bar{a}_3\bar{a}_4\wedge 1_M)(1_{E_5}\wedge \varepsilon\bar{r})d(f_2ij),\quad\hbox{by (3.7)} \\
 = &0.
\end{array}
$$
Hence, by (1.2) we have
$$
(1_{E_3}\wedge u_6)f_3^{\prime}=(1_{E_3}\wedge j^{\prime})f_4^{\prime}
\eqno{(3.28)}
$$
with $f_4^{\prime}\in [\Sigma^{p^nq+(p+2)q+1}M,E_3\wedge K]$.

Similarly, by Lemma 3.3 we have $\varepsilon(1_{K^{\prime}}\wedge
i)\upsilon j^{\prime}=-2j^{\prime}\alpha^{\prime}$. Then we have
$$\begin{array}{ll}
   & -2(1_{E_3}\wedge j^{\prime}\alpha^{\prime})f_4^{\prime} \\
= & (1_{E_3}\wedge \varepsilon(1_{K^{\prime}}\wedge i)
\upsilon j^{\prime})f_4^{\prime} \\
= & (1_{E_3}\wedge \varepsilon(1_{K^{\prime}}\wedge i)\upsilon)(1_{E_3}
\wedge j^{\prime})f_4^{\prime},\quad\hbox{by (3.28)} \\
 = &(1_{E_3}\wedge \varepsilon(1_{K^{\prime}}\wedge i)\upsilon)(1_{E_3}\wedge u_6)f_3^{\prime}
 ,\quad\hbox{by (3.24)}  \\
= & (1_{E_3}\wedge \varepsilon(1_{K^{\prime}}\wedge ij)
(\upsilon\wedge 1_M)\bar{m}_M)(1_{E_3}\wedge u_6)f_3^{\prime}
\\
= & (1_{E_3}\wedge \varepsilon(1_{K^{\prime}}\wedge ij))(1_{E_3}\wedge
\upsilon\wedge 1_M\bar{m}_M u_6)f_3^{\prime},\quad\hbox{by (3.24)} \\
=&(1_{E_3}\wedge \varepsilon(1_{K^{\prime}}\wedge ij))(1_{E_3}\wedge u_4 u_5)f_3^{\prime},
\quad\hbox{by (3.9)}    \\
= & (1_{E_3}\wedge \varepsilon(1_{K^{\prime}}\wedge ij)u_4 u_5)f_3^{\prime}
\\
= & 0.
\end{array}
$$
Thus, by (3.6) we have
$$
f_4^{\prime}=(1_{E_3}\wedge\pi)f_5^{\prime}
\eqno{(3.29)}
$$
with $f_5^{\prime}\in [\Sigma^{p^nq+(p+2)q+1}M, E_3\wedge
K^{\prime}\wedge M]$. This completes the proof of Lemma 3.11.

{\protect\vspace{-15pt}\rightline{$\square$}}

{\bf Lemma 3.12}\quad {\it For the above $f_5^{\prime}\in
[\Sigma^{p^nq+(p+2)q+1}M, E_3\wedge K^{\prime}\wedge M]$, we have $$
(\bar{b}_3\wedge 1_{K^{\prime}\wedge M})f_5^{\prime}=0.$$}

{\it Proof.}\quad The proof will be given later.

{\protect\vspace{-15pt}\rightline{$\square$}}

Now we give the proof of Theorem I.

{\bf Proof of Theorem I.} From Lemma 3.12, we have
$$
(\bar{b}_3\wedge 1_{K^{\prime}\wedge M})f_5^{\prime}=0.
\eqno{(3.30)}
$$
By virtue of (3.2), we have
$$
f_5^{\prime}=(\bar{a}_3\wedge 1_{K^{\prime}\wedge M})f_6^{\prime}
\eqno{(3.31)}
$$
with $f_6^{\prime}\in [\Sigma^{p^nq+(p+2)q+2}M, E_4\wedge
K^{\prime}\wedge M]$. By (3.27) and (3.2), we have
$$\begin{array}{ll}
  &(\bar{a}_2{\bar{a}}_3{\bar{a}}_4
  \wedge 1_{K^{\prime}\wedge M})(1_{E_5}\wedge\bar{r})d(f_2ij)
    \\
  = &(\bar{a}_2\wedge 1_{K^{\prime}\wedge M})(1_{E_3}
  \wedge (\upsilon\wedge 1_M){\bar{m}}_M u_6)f_3^{\prime}
  \\
= &(\bar{a}_2\wedge 1_{K^{\prime}\wedge M})(1_{E_3}
\wedge (\upsilon\wedge 1_M)\bar{m}_M)(1_{E_3}\wedge u_6)f_3^{\prime},\quad\hbox{by (3.28)}  \\
=&(\bar{a}_2\wedge 1_{K^{\prime}\wedge M})(1_{E_3}\wedge
(\upsilon\wedge 1_M)\bar{m}_M)
(1_{E_3}\wedge j^{\prime})f_4^{\prime},\quad\hbox{by (3.29)} \\
=&(\bar{a}_2\wedge 1_{K^{\prime}\wedge M})(1_{E_3}\wedge
(\upsilon\wedge 1_M)\bar{m}_M)(1_{E_3}\wedge j^{\prime})
(1_{E_3}\wedge\pi)f_5^{\prime},\quad\hbox{by (3.31)}\\
= & (\bar{a}_2\wedge 1_{K^{\prime}\wedge M})(1_{E_3}\wedge
(\upsilon\wedge 1_M)\bar{m}_M)(1_{E_3}\wedge j^{\prime})
(1_{E_3}\wedge\pi)(1_{E_3}\wedge\pi)(\bar{a}_3 \wedge
1_{K^{\prime}\wedge M})f_6^{\prime}
\\
=&(\bar{a_2}\bar{a}_3\wedge 1_{K^{\prime}\wedge
M})(1_{E_4}\wedge(\upsilon\wedge
1_M)\bar{m}_Mj^{\prime})(1_{E_4}\wedge \pi)f_6^{\prime}.
\end{array}
$$
That is,
$$
(\bar{a}_2{\bar{a}}_3{\bar{a}}_4\wedge 1_{K^{\prime}\wedge
M})(1_{E_5}\wedge\bar{r})d(f_2ij) =(\bar{a_2}\bar{a}_3\wedge
1_{K^{\prime}\wedge M})(1_{E_4}\wedge(\upsilon\wedge
1_M)\bar{m}_Mj^{\prime})(1_{E_4}\wedge \pi)f_6^{\prime}.
\eqno{(3.32)}
$$
Since $[(\bar{b}_4\wedge 1_K)(1_{E_4}\wedge \pi)f_6^{\prime}]\in
Ext_A^{4,p^nq+(p+2)q+2}(H^{\ast}K, H^{\ast}M)=0$ (cf. [7,
Proposition 2.2]), then by (3.1), we know that the $d_1$-cycle
$(\bar{b}_4\wedge 1_K)(1_{E_4}\wedge \pi)f_6^{\prime}$ is a
$d_1$-boundary. It follows that $(\bar{b}_4\wedge 1_K)(1_{E_4}\wedge
\pi)f_6^{\prime}=(\bar{b}_4\wedge 1_K)(\bar{c}_3\wedge
1_K)f_7^{\prime}$ for some $f_7^{\prime}\in
[\Sigma^{p^nq+(p+2)q+2}M, KG_3\wedge K]$. Thus we have
$$
(1_{E_4}\wedge \pi)f_6^{\prime}=(\bar{c}_3\wedge
1_K)f_7^{\prime}+(\bar{a}_4\wedge 1_K)f_8^{\prime}
\eqno{(3.33)}
$$
with $f_8^{\prime}\in[\Sigma^{p^nq+(p+2)q+3}M, E_5\wedge K]$.
Then by (3.32), (3.33) and (3.2), we have
$$
(\bar{a}_2\bar{a}_3\bar{a}_4\wedge 1_{K^{\prime}\wedge
M})(1_{E_5}\wedge\bar{r})d(f_2ij)=(\bar{a}_2\bar{a}_3\bar{a}_4\wedge
1_{K^{\prime}\wedge M}) (1_{E_5}\wedge (\upsilon\wedge
1_M)\bar{m}_Mj^{\prime})f_8^{\prime}. \eqno{(3.34)}
$$
Moreover, by composing $(1_{E_2}\wedge \bar{r})$ on (3.13) it is
easy to get that
$$
(\bar{a}_2\bar{a}_3\bar{a}_4\wedge 1_{K^{\prime}\wedge
M})(1_{E_5}\wedge\bar{r})d(f_2ij)
=(1_{E_2}\wedge \bar{r}(i^{\prime\prime}\wedge
1_K)\beta)d(\eta_{n,2}^{\prime}i^{\prime}ij).
\eqno{(3.35)}
$$
Combining (3.34) and (3.35) yields
$$
(1_{E_2}\wedge \bar{r}(i^{\prime\prime}\wedge
1_K)\beta)d(\eta_{n,2}^{\prime}i^{\prime}ij)
=(\bar{a}_2\bar{a}_3\bar{a}_4\wedge 1_{K^{\prime} \wedge
M})(1_{E_5}\wedge (\upsilon\wedge
1_M)\bar{m}_Mj^{\prime})f_8^{\prime}. \eqno{(3.36)}
$$
Notice that $\bar{r}(i^{\prime\prime}\wedge 1_K)=(\upsilon\wedge
1_M)\bar{m}_Mj^{\prime}$ (see Lemma 3.3). Then (3.36) can turn into
$$
(1_{E_2}\wedge (\upsilon\wedge
1_M)\bar{m}_Mj^{\prime}\beta)d(\eta_{n,2}^{\prime}i^{\prime}ij)
=(1_{E_2}\wedge (\upsilon\wedge
1_M)\bar{m}_Mj^{\prime})(\bar{a}_2\bar{a}_3\bar{a}_4\wedge
1_K)f_8^{\prime}. \eqno{(3.37)}
$$
By (3.37) and (3.6), we have
$$
(1_{E_2}\wedge j^{\prime}\beta)d(\eta_{n,2}^{\prime}i^{\prime}ij)
=(1_{E_2}\wedge j^{\prime})(\bar{a}_2\bar{a}_3\bar{a}_4\wedge
1_K)f_8^{\prime}+(1_{E_2}\wedge
j^{\prime}\alpha^{\prime})f_9^{\prime}
\eqno{(3.38)}
$$
with $f_9^{\prime}\in [\Sigma^{p^nq+(p+1)q+1}M,E_2\wedge K]$. From
[7, p.489], we know that the left hand side of (3.38) has filtration
4. However, since the first term of the right hand side of (3.38)
has filtration $\geq 5$, then the second term of (3.38) must be of
filtration $4$. So $f_9^{\prime}$ has filtration $\leq 3$. Notice
the facts that $Ext_A^{3,p^nq+(p+1)q+2}(H^{\ast}K,H^{\ast}M)=0$ (cf.
Proposition 2.5) and $Ext_A^{2,p^nq+(p+1)q+1}(H^{\ast}K,H^{\ast}M)
\cong \mathbb{Z}_p\{\beta_{\ast} i^{\prime}_{\ast}(\tilde{h}_n)\}$
(cf. Proposition 2.7), then we have $(\bar{b}_2\wedge
1_K)f_9^{\prime} =(1_{KG_2}\wedge \beta)(1_{KG_2}\wedge
i^{\prime})(\tilde{h}_n)$. Let $\xi_n=(\bar{a}_0\bar{a}_1\wedge
1_K)f_9^{\prime}$. Then $\xi_n$ is represented by
$\beta_{\ast}i_{\ast}^{\prime}(\tilde{h}_n)$ in the Adams spectral
sequence. And so $\zeta_n=\xi_ni$ is represented by
$i^{\ast}\beta_{\ast}i_{\ast}^{\prime}(\tilde{h}_n)=
\beta_{\ast}i_{\ast}^{\prime}i^{\ast}(\tilde{h}_n)=\beta_{\ast}i_{\ast}^{\prime}i_{\ast}(h_n)
\not=0\in Ext_A^{2,p^nq+(p+1)q+1}(H^{\ast}K,\mathbb{Z}_p)$ {(cf.
Proposition 2.6)}. Thus Theorem I is proved.

{\bf Proof of Lemma 3.12.}\quad We first recall three cofibrations
given in [7]
$$
\Sigma^{-1}K\stackrel{\upsilon
  j^{\prime}}{\longrightarrow}\Sigma^{q}K^{\prime}\stackrel{\bar{\psi}}{\longrightarrow}
  K_2^{\prime}\stackrel{\bar{\rho}}{\longrightarrow}K\quad\hbox{(see
  [7, (2.5)])},
  \eqno{(3.39)}$$
$$
\Sigma^{q-1}K_1^{\prime}\stackrel{\varepsilon(1_{K_1^{\prime}}\wedge
i)}{\longrightarrow}M\stackrel{w_2}{\longrightarrow}
X\stackrel{u_2}{\longrightarrow}\Sigma^{q}K^{\prime}\quad\hbox{(see
[7, (3.7)])} , \eqno{(3.40)}
$$
$$X\stackrel{\bar{\psi}u_2}{\longrightarrow}K_2^{\prime}\stackrel{w_3}{\longrightarrow}
K^{\prime}\wedge W\stackrel{u_3}{\longrightarrow}\Sigma
X\quad\hbox{(see [7, (3.10)])} \eqno{(3.41)}
$$
with the relation $u_2u_3=-\upsilon j^{\prime}\pi$[7, (3.11)]. By
composing $(\bar{a}_2\wedge 1_{K^{\prime}\wedge M})$ on (3.27), we
have
$$\begin{array}{ll}
 &(\bar{a}_2\bar{a}_3\bar{a}_4\wedge 1_{K^{\prime}\wedge M})(1_{E_5}\wedge \bar{r})d(f_2ij) \\
=&(\bar{a}_2\wedge 1_{K^{\prime}\wedge
M})(1_{E_3}\wedge(\upsilon\wedge
1_M)\bar{m}_Mu_6)f_3^{\prime},\quad\hbox{by (3.28)}\\
=&(\bar{a}_2\wedge 1_{K^{\prime}\wedge
M})(1_{E_3}\wedge(\upsilon\wedge 1_M)\bar{m}_M)(1_{E_3}
\wedge j^{\prime})f_4^{\prime},\quad\hbox{by (3.29)}\\
=&(\bar{a}_2\wedge 1_{K^{\prime}\wedge
M})(1_{E_3}\wedge(\upsilon\wedge 1_M)\bar{m}_M)(1_{E_3}\wedge
j^{\prime})(1_{E_3}\wedge \pi)f_5^{\prime}.
\end{array}$$
That is,
$$
(\bar{a}_2\bar{a}_3\bar{a}_4\wedge 1_{K^{\prime}\wedge
M})(1_{E_5}\wedge \bar{r})d(f_2ij) =(\bar{a}_2\wedge
1_{K^{\prime}\wedge M})(1_{E_3}\wedge(\upsilon\wedge
1_M)\bar{m}_M)(1_{E_3}\wedge j^{\prime})(1_{E_3}\wedge
\pi)f_5^{\prime}. \eqno{(3.42)}$$ By composing $(1_{E_2}\wedge
(1_{K^{\prime}}\wedge j))$ on (3.42), we have
$$
\begin{array}{cl}
&(1_{E_2}\wedge (1_{K^{\prime}}\wedge j))
(\bar{a}_2\bar{a}_3\bar{a}_4\wedge 1_{K^{\prime}\wedge M})(1_{E_5}\wedge
\bar{r})d(f_2ij)\\
=&(1_{E_2}\wedge (1_{K^{\prime}}\wedge j))
 (\bar{a}_2\wedge 1_{K^{\prime}\wedge M})(1_{E_3}\wedge(\upsilon\wedge
 1_M)\bar{m}_M)
 (1_{E_3}\wedge
 j^{\prime})(1_{E_3}\wedge \pi)f_5^{\prime}.
 \end{array}
 \eqno{(3.43)}
 $$
On the one hand, for the left hand side of (3.43), we have
$$\begin{array}{ll}
  &(1_{E_2}\wedge (1_{K^{\prime}}\wedge j))
   (\bar{a}_2\bar{a}_3\bar{a}_4\wedge 1_{K^{\prime}\wedge M})(1_{E_5}\wedge \bar{r})d(f_2ij) \\
  = &-(\bar{a}_2\bar{a}_3\bar{a}_4
  \wedge 1_{K^{\prime}})(1_{E_5}\wedge(1_{K^{\prime}}\wedge j)\bar{r})d(f_2ij),\quad\hbox{by
  Lemma 3.4}  \\
  = &-(\bar{a}_2\bar{a}_3\bar{a}_4
  \wedge 1_{K^{\prime}})(1_{E_5}\wedge \bar{\Delta}_{K^{\prime}})d(f_2ij)
   .
\end{array}$$
On the other hand, for the right hand side of (3.43) we have
$$\begin{array}{ll}
   &(1_{E_2}\wedge (1_{K^{\prime}}\wedge j))
    (\bar{a}_2\wedge 1_{K^{\prime}\wedge M})(1_{E_3}\wedge(\upsilon\wedge
    1_M)\bar{m}_M)
    (1_{E_3}\wedge
    j^{\prime})(1_{E_3}\wedge \pi)f_5^{\prime}.
    \\
  = & -(\bar{a}_2\wedge 1_{K^{\prime}})(1_{E_3}
\wedge(1_{K^{\prime}}\wedge j)(\upsilon\wedge 1_M)\bar{m}_Mj^{\prime}\pi)f_5^{\prime} \\
= & -(\bar{a}_2\wedge 1_{K^{\prime}})(1_{E_3}\wedge (\upsilon\wedge 1_{s^0})(1_M\wedge j)
\bar{m}_Mj^{\prime}\pi)f_5^{\prime},\quad\hbox{since $(1_M\wedge j)\bar{m}_M=1_M$} \\
= & -(\bar{a}_2\wedge 1_{K^{\prime}})(1_{E_3}\wedge \upsilon j^{\prime}\pi)f_5^{\prime},
\quad\hbox{since $u_2u_3=-\upsilon j^{\prime}\pi$} \\
= &(\bar{a}_2\wedge 1_{K^{\prime}})(1_{E_3}\wedge
u_2u_3)f_5^{\prime}.
\end{array}$$
Thus we have
$$
(\bar{a}_2\bar{a}_3\bar{a}_4
\wedge 1_{K^{\prime}})(1_{E_5}\wedge \bar{\Delta}_{K^{\prime}})d(f_2ij)
=-(\bar{a}_2\wedge 1_{K^{\prime}})(1_{E_3}\wedge
u_2u_3)f_5^{\prime}.
\eqno{(3.44)}
$$

Let $X$ be the cofibre of
$\varepsilon(1_{K^{\prime}}\wedge i):\Sigma^{q-1}K^{\prime}\longrightarrow M$
given by the cofibration
$$\Sigma^{q-1}K^{\prime}\stackrel{\varepsilon(1_{K^{\prime}}\wedge i)}{\longrightarrow}
M\stackrel{w_2}{\longrightarrow}
X\stackrel{u_2}{\longrightarrow}\Sigma^q K^{\prime}\quad\hbox{(see
[7, (3.7)])}.$$

[7,(3.8)] follows from [7, (3.7)] and [7, (3.6)] that
$$(\bar{a}_4\wedge 1_{K^{\prime}})(1_{E_5}\wedge \bar{\Delta}_{K^{\prime}})
d(f_2ij)=(1_{E_4}\wedge u_2)f_3$$ for some $f_3\in
[\Sigma^{p^nq+(p+2)q+1}M, E_4 \wedge X]$. By composing
$(\bar{a}_2\bar{a}_3\wedge 1_{K^{\prime}})$ on [7, (3.8)], we have
$$
(\bar{a}_2\bar{a}_3\bar{a}_4\wedge 1_{K^{\prime}})(1_{E_5}\wedge
\bar{\Delta}_{K^{\prime}})d(f_2ij) =(\bar{a}_2\bar{a}_3\wedge
1_{K^{\prime}})(1_{E_4}\wedge u_2)f_3. \eqno{(3.45)}
$$
 Combining
(3.44) and (3.45) yields
$$
(\bar{a}_2\wedge 1_{K^{\prime}})(1_{E_3}\wedge
u_2u_3)f_5^{\prime}=-(\bar{a}_2\bar{a}_3\wedge 1_{K^{\prime}})(1_{E_4}\wedge
u_2)f_3.
\eqno{(3.46)}
$$
By [10, (1.7)], (3.46) can turn into
$$
(1_{E_2}\wedge u_2)(\bar{a}_2\wedge 1_X)(1_{E_3}\wedge u_3)f_5^{\prime}=-(1_{E_2}\wedge
u_2)(\bar{a}_2\bar{a}_3\wedge 1_X)f_3.
\eqno{(3.47)}
$$
From (3.47) and (3.40) we have
$$
(\bar{a}_2\wedge 1_X)(1_{E_3}\wedge u_3)f_5^{\prime}=-(\bar{a}_2\bar{a}_3\wedge 1_X)f_3
+(1_{E_2}\wedge w_2)\bar{f}_4
\eqno{(3.48)}
$$
with $\bar{f}_4\in [\Sigma^{p^nq+(p+2)q-1}M, E_2\wedge M]$.

Note
that $(\bar{b}_2\wedge 1_M)\bar{f}_4\in [\Sigma^{p^nq+(p+2)q-1}M, KG_2\wedge
M]=0$ by the exact sequence
$[\Sigma^{p^nq+(p+2)q-1}M, KG_2]\stackrel{(1\wedge i)_{\ast}}{\longrightarrow}
[\Sigma^{p^nq+(p+2)q-1}M, KG_2\wedge M]\stackrel{(1\wedge j)_{\ast}}{\longrightarrow}
[\Sigma^{p^nq+(p+2)q-2}M,KG_2]$ induced by (1.1), where the first
and the last group are zero by the fact that $\pi_{p^nq+(p+2)q+r}KG_2\cong
Ext_A^{2,p^nq+(p+2)q+r}(\mathbb{Z}_p,\mathbb{Z}_p)=0
$ for $r=0,-1,-2$ (cf. [2]). Hence, $\bar{f}_4=(\bar{a}_2\wedge
1_M)\bar{f}_5$ for some $\bar{f}_5\in [\Sigma^{p^nq+(p+2)q}M,E_3\wedge
M]$. By (3.2) and (3.48), we have
$$
(1_{E_3}\wedge u_3)f_5^{\prime}=
-(\bar{a}_3\wedge 1_X)f_3+(1_{E_3}\wedge w_2)\bar{f}_5^{\prime}+(\bar{c}_2\wedge
1_{X})g_6
\eqno{(3.49)}$$
with $g_6\in [\Sigma^{p^nq+(p+2)q}M, KG_2\wedge X]$. And so
we have
$$
(\bar{b}_3\wedge 1_X)(1_{E_3}\wedge
u_3)f_5^{\prime}=(\bar{b}_3\wedge 1_X)(1_{E_3}\wedge
w_2)\bar{f}_5+(\bar{b}_3\bar{c}_2\wedge 1_X)g_6.\eqno{(3.50)}$$
Since $Ext_A^{3,p^nq+(p+2)q}(H^{\ast}M,H^{\ast}M) \cong
\mathbb{Z}_p\{i_{\ast}j_{\ast}\overline{\overline{h_ng_0}},
j^{\ast}i^{\ast}\overline{\overline{h_ng_0}}\}$ (cf. Proposition
2.3), then we can have $(\bar{b}_3\wedge
1_M)\bar{f}_5=\lambda_1\overline{\overline{h_ng_0}}ij+
\lambda_2(1_{KG_3}\wedge ij)\overline{\overline{h_ng_0}}$ for some
$\lambda_1,\lambda_2\in \mathbb{Z}_p$, where
$\overline{\overline{h_ng_0}} \in
[\Sigma^{p^nq+(p+2)q+1}M,KG_3\wedge M]$. And so
$$0=\lambda_1(\bar{c}_3\wedge 1_M)\overline{\overline{h_ng_0}}ij+
\lambda_2(\bar{c}_3\wedge 1_M)(1_{KG_3}\wedge
ij)\overline{\overline{h_ng_0}}.$$
By composing $i$ on the
 above equality,
we get that $\lambda_2(\bar{c}_3\wedge 1_M)(1_{KG_3}\wedge
ij)\overline{\overline{h_ng_0}}i=0$. However, since
$d_2(i^{\ast}(ij)_{\ast}\overline{\overline{h_ng_0}})
=i^{\ast}d_2(i_{\ast}j_{\ast}\overline{\overline{h_ng_0}})\not=0$
(see Proposition 2.4 (1)), we get that $(\bar{c}_3\wedge
1_M)(1_{KG_3}\wedge ij)\overline{\overline{h_ng_0}}i \not=0$. Thus,
we have that $$\lambda_2=0, \lambda_1(\bar{c}_3\wedge 1_M)
\overline{\overline{h_ng_0}}ij=0.$$ Note that
$d_2(j^{\ast}i^{\ast}\overline{\overline{h_ng_0}})\not=0$ by
Proposition 2.4 (2), then $(\bar{c}_3\wedge
1_M)\overline{\overline{h_ng_0}}ij\not=0$. Thus we see that
$\lambda_1=0$. From the above discussion, we know that
$(\bar{b}_3\wedge 1_M)\bar{f}_5=0$ and (3.50) can turn into
$$
(\bar{b}_3\wedge 1_X)(1_{E_3}\wedge u_3)f_5^{\prime}=(\bar{b}_3\bar{c}_2\wedge
1_X)g_6.
\eqno{(3.51)}
$$
The argument of the proof from [7, (3.16)] to [7, p.491] shows that
$(\bar{b}_3\wedge 1_X)(1_{E_3}\wedge u_3)f_6=-(\bar{b}_3\bar{c}_2
\wedge 1_X)\tilde{l}_0$ in [7, (3.16)] implies $(\bar{b}_3\wedge
1_{K^{\prime}\wedge M})f_6=0$. By a similar argument as in [7], we
can also show that (3.51) implies (3.30) holds.

{\protect\vspace{-15pt}\rightline{$\square$}}

{\it\bf Proof of Theorem II}\quad By Theorem I , we get
that
$$\beta_{\ast}{i^{\prime}}_{\ast}i_{\ast}({h}_n)\not=0\in
Ext_A^{2,p^nq+(p+1)q+1}(H^{\ast}K,\mathbb{Z}_p)$$
is a permanent cycle in the Adams spectral sequence and converges to a
nontrivial element $\zeta_n \in \pi_{p^nq+(p+1)q-1}K$.

Consider
the following composition of maps
$$\bar{f}:\Sigma^{p^nq+(p+1)q-1}S\stackrel{\zeta_n}{\longrightarrow}K
\stackrel{jj^{\prime}\beta}{\longrightarrow}
\Sigma^{-pq+2}S.$$
Since $\zeta_n$ is represented up to nonzero scalar
by $\beta_{\ast}{i^{\prime}}_{\ast}i_{\ast}({h}_n)\in
Ext_A^{2,p^nq+(p+1)q+1}(H^{\ast}K,\mathbb{Z}_p)$ in the Adams spectral sequence,
then the above $\bar{f}$ is
represented up to nonzero scalar by
$
\bar{c}=(jj^{\prime}\beta)_{\ast}\beta_{\ast}{i^{\prime}}_{\ast}i_{\ast}({h}_n)
$ in the Adams spectral sequence.

Meanwhile, it is well known that the $\beta$-element
$\beta_2=jj^{\prime}\beta^2 i^{\prime}i$ is represented by $k_0\in
Ext_A^{2,2pq+q} (\mathbb{Z}_p,\mathbb{Z}_p)$ in the Adams spectral
sequence. By the knowledge of Yoneda products we can see that
$\bar{f}$ is represented (up to nonzero scalar) by $$\bar{c}
=k_0h_n\not= 0\in Ext_A^{3,q(p^n+2p+1)}(\mathbb{Z}_p,\mathbb{Z}_p)$$
in the Adams spectral sequence(cf. [6, Table 8.1]).

Moreover, we know that
$Ext_A^{3-r,q(p^n+2p+1))+(-r+1)}(\mathbb{Z}_p,\mathbb{Z}_p)=0$ for
$r\geq 2$, then $k_0h_n$ cannot be hit by any differential in the
Adams spectral sequence, and so the corresponding homotopy element
$\bar{f}\in\pi_{\ast}S$ is nontrivial and of order $p$. This
finishes the proof of Theorem II.

{\protect\vspace{-15pt}\rightline{$\square$}}

\begin{center}
{\bf\large References}
\end{center}
\begin{enumerate}
\item Cohen R. Odd primary families in stable homotopy theory.
{\it Memoirs of the American Mathematical Society}, 1981, {\bf 242}.
\item Liulevicius A. The factorizations of cyclic reduced powers by
secondary
cohomology operations. {\it Memoirs of the American Mathematical Society}, 1962, {\bf 42}.
\item Wang X J. Ordered cochain complex and $H^{4,\ast}(A_p)$ for
$p\geq 5$. {\it Beijing Mathematics}, 1995, 1: $80\sim 99$.
\item Liu X G. A nontrivial product in the stable homotopy groups of spheres.
{\it Sci. in China Ser. A}, 2004, {\bf 47}(6): $831\sim 841$.
\item H. Toda. On spectra realizing exterior part of the Steenrod algebra.
{\it Topology}, 1971, {\bf 10}: $53\sim 65$.
\item Aikawa T. 3-dimensional cohomology of the mod p Steenrod
algebra. {\it Mathematica Scandinavica}, 1980, {\bf 47}: $91\sim 115$.
\item Lin J K. A new family of filtration three in the stable
homotopy
of spheres. {\it Hiroshima Mathematical Journal}, 2001, {\bf 31}: $477\sim 492$.
\item Lin J K and Zheng Q B. A new family of filtration seven in the
stable homotopy
of spheres. {\it Hiroshima Mathematical Journal}, 1998, {\bf 28}: $183\sim 205$.
\item Cohen R. and Goerss P. Secondary cohomology operations that
detect
homotopy classes,
{\it Topology}, 1984, {\bf 23}: $177-194$.
\item Toda H. Algebra of stable homotopy of $\mathbb{Z}_p$-spaces and
applications. {\it
Journal of Mathematics of Kyoto University}, 1971, {\bf 11}:
$197\sim 251$.
\item  Oka S. Multilicative structure of finite ring spectra and
stable homotopy of
spheres. {\it Algebraic Topology (Aarhus), Lecture
Notes in Mathematics}, 1984, {\bf 1051}: $418\sim 441$.
Springer-Verlag.
\item Thomas E. and Zahler R. Generalized higher order cohomology
operations
 and stable homotopy groups of spheres. {\it Advances in Mathematics}, 1976, {\bf
20}:
$287\sim 328$.
\end{enumerate}
\end{document}